\documentclass[10pt,a4paper,article]{amsart}
\usepackage{amssymb,stmaryrd}
\usepackage{amsthm,amstext}
\usepackage{amssymb,amsmath}
\usepackage{amsfonts}
\usepackage{MnSymbol,slashed}
\usepackage{fancyhdr}
\usepackage{graphicx}
\usepackage{pdfpages}
\usepackage{float}
\usepackage{subfig}
\usepackage{simplewick}
\usepackage{graphics,psfrag}
\usepackage{pstricks}
\usepackage{latexsym}
\usepackage{changepage}
\usepackage{caption}
\usepackage{paralist}
\allowdisplaybreaks[1]

\theoremstyle{plain}
\newtheorem{theorem}{Theorem}[section] 

\newtheorem{lemma}[theorem]{Lemma}
\newtheorem{corollary}[theorem]{Corollary}
\newtheorem{proposition}[theorem]{Proposition}
\newtheorem{definition}[theorem]{Definition}
\newtheorem{example}[theorem]{Example}
\newtheorem{remark}[theorem]{Remark}
\newcommand{\CC}{\hbox{{$\mathcal C$}}}

\newcommand{\CF}{\hbox{{$\mathcal F$}}}

\newcommand{\CS}{{\hbox{{$\mathcal S$}}}}

\newcommand{\CT}{\hbox{{$\mathcal T$}}}

\newcommand{\CR}{\hbox{{$\mathcal R$}}}
\newcommand{\CQ}{\hbox{{$\mathcal Q$}}}

\newcommand{\C}{\mathbb{C}}
\newcommand{\R}{\mathbb{R}}
\newcommand{\F}{\mathbb{F}}
\newcommand{\Z}{\mathbb{Z}}

\newcommand{\cg}{\hbox{{$\mathfrak g$}}}
\newcommand{\ch}{\hbox{{$\mathfrak h$}}}

\newcommand{\isom}{{\cong}}
\newcommand{\eps}{{\epsilon}}
\newcommand{\tens}{\mathop{{\otimes}}}
\newcommand{\la}{{\triangleright}}

\newcommand{\id}{\mathrm{id}}

\renewcommand{\>}{\rangle}

\newcommand{\prodf}{\includegraphics{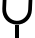}}
\newcommand{\coprodf}{\includegraphics{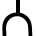}}
\allowdisplaybreaks[1]

\begin{document}

\author{Shahn Majid}
\address{School of Mathematical Sciences\\ Queen Mary University of London \\ Mile End Rd, London E1 4NS \\ 
and Cambridge Quantum Computing}
\email{ s.majid@qmul.ac.uk}
\thanks{Work done on sabbatical at Cambridge Quantum Computing}

\title[Quantum and braided ZX calculus]{Quantum and braided ZX calculus}
	\begin{abstract} We revisit the notion of interacting Frobenius Hopf algebras for ZX-calculus in quantum computing, with focus on allowing the algebras to be noncommutative and coalgebras to be noncocommutative.  We  introduce the notion of $*$-structures in ZX-calculus at this algebraic level and construct examples based on the quantum group $u_q(sl_2)$ at a root of unity. We provide an abstract formulation of the Hadamard gate related to Hopf algebra self-duality. We then solve the problem of extending the notion of interacting Hopf algebras and ZX-calculus to take place in a braided tensor category. In the ribbon case, the Hadamard gate coming from braided self-duality obeys a modular identity. We give the example of $b_q(sl_2)$, the self-dual braided version of $u_q(sl_2)$. \end{abstract}
\maketitle

\section{Introduction}

This work is on the interface between Computer Science in the area of quantum computing, which is based particularly on a diagrammatic formulation of ZX calculus\cite{CD}, and the theory of Hopf algebras or quantum groups as understood in algebra and mathematical physics. The preliminaries in Section~\ref{secquaZX} start with a concise but self-contained reworking of the central algebraic structure, namely that of a Frobenius Hopf algebra or `interacting pair of Hopf algebras', but without assuming that  our algebras are commutative or coalgebras cocommutative or our Frobenius forms symmetric. Here, Corollary~\ref{corHH*} revisits a recent result in \cite{CoD} that every finite-dimensional Hopf algebra $H$ `amplifies' to such a pair, but with a more direct treatment that extends to the braided case and identifies the other member of the pair as $H^{*op}$. We keep track of braid crossings in the diagrammatic proofs as preparation for the braided theory later, and we compute details for the example of $u_q(sl_2)$ at $q$ a primitive $n$th root of unity. The simplest case is $n=2$ which, in our conventions from \cite{AziMa}, is noncommutative and noncocommutative even for $q=-1$ and given in detail. 

Section~\ref{secstar} introduces $*$-structures on Frobenius and Frobenius-Hopf algebras. The nicest assumptions are where $H$ but require unimodularity. For $u_q(sl_2)$ we show how this gets modified, with its natural $*$ operation making it a {\rm flip}-$*$ Hopf algebra. Phases and an abstract formulation of a `Hadamard gate' are also discussed, the latter in three types. Here, Type 1 is the more obvious one but Type 2 is equivalent to $H$ self-dual as a Hopf algebra and Type 3 to $H$ anti-self dual, which are therefore more natural from a Hopf algebra point of view. $u_{-1}(sl_2)$ and the reduced Taft algebra $u_q(b_+)\subset u_q(sl_2)$ provide examples. 

The main new result of the paper is in Section~\ref{secbraZX}, where we construct interacting Hopf algebras and ZX calculus in a braided tensor category. The general result in the braided case that every braided Hopf algebra with suitable integrals amplifies to a braided Frobenius-Hopf algebra (Corollary~\ref{corbraHH*}) is read off from the diagrammatic proofs in Section~\ref{secquaZX}, so that the main additional work is the construction of examples. We do this (Proposition~\ref{proptransH}) by a process of transmutation~\cite{Ma:bra,Ma:tra} whereby every quasitriangular Hopf algebra in the sense of Drinfeld\cite{Dri,Ma} has a braided version. We include detailed computations for $b_q(sl_2)$ obtained by transmutation from $u_q(sl_2)$. Our Type 2 `braided hadamard gate' $\ch$ here is essentially the braided Fourier transform for this class of braided Hopf algebras,  shown in \cite{LyuMa} along with the ribbon structure to obey the modular identities needed for topological invariants. 

We conclude this introduction with a little of the background from Computer Science. In quantum computing, any unitary operator can be regarded as a `quantum gate', but the useful content is to identify a specific collection of such `gates' from which others can be built or approximated to arbitrary accuracy. To this end, let $V=\C^2$ for a singe `qubit' system and consider unitaries on the $m$-qubit system $V^{\tens m}$. Then some useful gates are the Hadamard gate 
\[ \ch: V\to V,\quad \ch={1\over\sqrt 2}\begin{pmatrix}1 & 1\\ 1 & -1\end{pmatrix} \]
and two flavours of maps  $\Delta_{\color{green}\bullet}:V\tens V\to V, \mu_{\color{red}\bullet}:V\to V\tens V$ (traditionally depicted with red and green nodes) which, up to normalisation, come from the  product and coproduct of the group Hopf algebra $\C\Z_2$ and its dual Hopf algebra\cite{CD}. One can further include versions of the nodes modified by phases. The merit of this particular collection of `gates' is that the Hopf algebra structures allow one to simplify and rewrite compositions according to certain identities, leading to an efficient framework for working with them. Moreover, composition of these linear maps can conveniently be expressed as `string diagrams' where the copies of $V$ are written side by side and maps applied as for a flow chart or `quantum circuit' flowing down the page. In this case, the core identities between the `gates' appear more abstractly at the diagram level as a graphical ZX-calculus\cite{CD} and the role of $\C \Z_2$ is to realise this concretely.  The role of $\C\Z_2$ and its dual can then be played by another `interacting pair' of Hopf algebras\cite{DD} namely $H$ and an associated Hopf algebra $\tilde H$ with green product $\mu_{\color{green}\bullet}$ and red coproduct $\Delta_{\color{red}\bullet}$. 

In these string diagrams as used in ZX calculus, there is no significance to the over or under crossings of strings. On the other hand, the idea of doing algebraic operations at a string diagram level has already been familiar in algebra since the early 1990s as the theory of `braided groups' or Hopf algebras in braided categories \cite{Ma:bra,Ma:bg,Ma:tra,Ma:cro,Ma:alg,Ma}. Here the `wiring diagram' approach to algebraic identities was essential to keep track of under- and over- crossings, as these correspond to typically different operations. In the simplest `anyonic' example, the braiding is a power a parameter $q$ and the inverse crossing at the same point would be require $q^{-1}$. The theory was needed as a way to understand deeper aspects of $q$-deformation quantum groups\cite{Dri}. This body of work  therefore begs the question as to whether ZX-calculus has a useful version at the level of braided tensor categories. Here, we show that at least for the core algebraic structures, namely the notion of a pair interacting Hopf algebras and Hadamard gate, this is indeed the case. The paper concludes with some remarks about further directions. 

\subsection*{Acknowledgements}  I thank the team at CQC for helpful discussions and references, notably B. Coecke, A. Cowtan, R. Duncan, A. Kissinger, Q. Wang. I also thank T. Fritz for a correction. 

\section{Construction of noncommutative interacting Hopf algebras}\label{secquaZX}

This section is based heavily on \cite{CD,DD,CoD} and does not claim to be particularly new, but we  are careful to ensure that where possible we do not assume  that our (co)algebras $A$ are (co)commutative or that our bilinear forms are symmetric. Section~\ref{secpre}, in particular, is strictly preliminary and included only to be self-contained for readers from a more algebraic background. Section~\ref{secFHopf} is a reworking of \cite{CoD} but we go further with new more direct proofs and  in which we keep track of braid crossings as preparation for Section~\ref{secbraZX}. We also strip out the formal categorical setting and focus on the content at an algebraic level. 

We will write our maps composed flowing downward and use juxtaposition of objects to denote $\tens$. Braidings are represented by braid crossings\cite{JoyStr} through which other morphisms can be pulled as part of the functoriality of the braiding. An introduction to doing actual algebra and Hopf algebra at this level can be found in \cite{Ma:alg,Ma,Ma:pri}.  We will give main formulae in both conventional and diagrammatic form. However, we will denote coproducts by $\Delta$ as usual in Hopf algebra theory, not by $\delta$ as favoured in the Computer Science literature. Adjunction of morphisms is done by rotating anticlockwise in the plane, which it seems is the opposite of the more common convention in Computer Science. Purely diagrammatic proofs work in any braided category, but for examples built from Hopf algebras we will assume we are in the $k$-linear case, where $k$ is a field.  

\subsection{Preliminaries on F-algebras}\label{secpre}

We recall that a unital algebra $A$ over a field $k$ is {\em Frobenius} if there is a bilinear form $(\ ,\ ):A\tens A\to k$ which is nondegenerate and such that $(ab,c)=(a,bc)$ for all $a,b,c\in A$. The latter says that the Frobenius bilinear form descends to a map $A\tens_A A\to k$ or equivalently a map $A\to k$. Nondegenerate as a linear map has its usually meaning which we write explicitly as existence of a {\em metric} $g=g^1\tens g^2\in A\tens A$ (sum of terms understood) such that  for all $a\in A$, 
\[ (a,g^1)g^2=a=g^1(g^2,a),\quad  \includegraphics[scale=0.6]{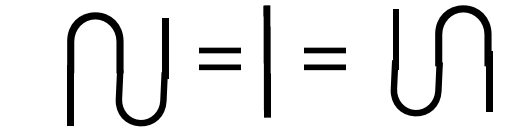}\]
where we also write the diagrammatic form with $(\ ,\ )=\cup$ and $g=\cap$. This makes $A$ a rigid object in the category of vector spaces and requires that $A$ is finite-dimensional. 

\begin{lemma}\label{gcentral} In a Frobenius algebra, $g$ commutes with the algebra product.\end{lemma}
\proof $g^1(g^2a,b)=g^1(g^2,ab)=ab=a g^1(g^2,b)$ for all $a,b\in A$, using the Frobenius and inverse properties. This implies $g^1(g^2a,g^{1'})\tens g^{2'}=ag^1(g^2,g^{1'})\tens g^{2'}$ or $g^1\tens (g^2a,g^{1'})g^{2'}=ag^1\tens (g^2,g^{1'}) g^{2'}$ and hence $ g^1\tens g^2 a=a g^1 \tens g^2$ for all $a\in A$. \endproof

In particular, the element $\mu(g)$, where we apply the product $\mu:A\tens A\to A$ of $A$ to $g$, is in the centre, $\mu(g)\in Z(A)$. Next, recall that reversing the axioms of a unital algebra gives the notion of a counital coalgebra $(A,\Delta,\eps)$ where the coproduct $\Delta:A\to A\tens A$ is coassociative and the counit $\eps:A\to k$ obeys $(\eps\tens\id)\Delta=\id=(\id\tens\eps)\Delta$. These have a standard diagrammatic representation in terms of tree diagrams and pruning identities\cite{Ma:alg}. The parallel notion  in Computer Science to a Frobenius algebra is at first sight different and we will therefore denote it differently, as an F-algebra.

\begin{definition}cf\cite{CD}\label{defFalg} An F-algebra is  a pair consisting of a unital algebra $(A,\mu,1)$ and counital coalgebra $(A,\Delta,\eps)$ where $\mu=\prodf$ and $\Delta=\coprodf$ are such that
\[ (\id\tens\mu)(\Delta\tens\id)=\Delta\mu=(\mu\tens\id)(\id\tens\Delta),\quad \includegraphics[scale=0.5]{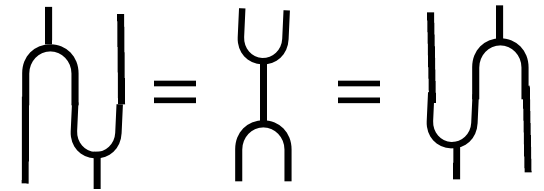}\]
An $F$-algebra is called {\em special} if $\mu\Delta=\id$ or \includegraphics[scale=0.3]{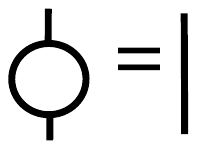} (or quasispecial if a nonzero multiple of this).
\end{definition}

\begin{figure}
 \[ \includegraphics[scale=1]{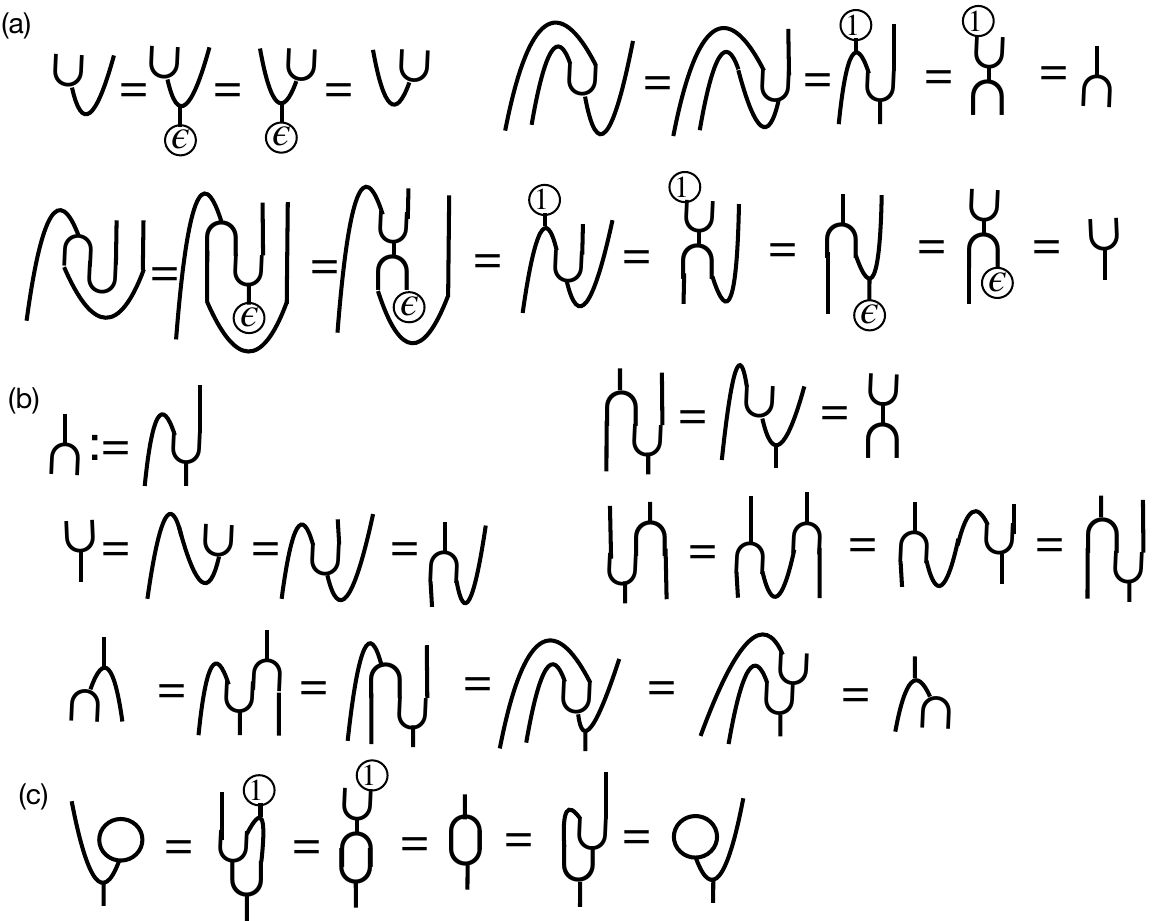}\]
 \caption{\label{figFalg} (a)-(b) Proof that F-algebras and Frobenius algebras are equivalent and (c) the special property as $\mu(g)=1$.}
 \end{figure}

\begin{lemma}\label{lemFalg}The following are equivalent for a unital algebra $A$

(1) A is the algebra part of an F-algebra.

(2) A is a Frobenius algebra.

\noindent The F-algebra is special if and only if $\mu(g)=1$. 
\end{lemma}
\proof (i) If we are given an F-algebra then let $(\ ,\ )=\cup=\eps\mu$ and $g=\cap=\Delta(1)$. The first half of the F-algebra axiom with $1$ on the top left leg and $\eps$ on the bottom right leg, recovers the first half of the axiom for a nondegenerate bilinear pairing. Similarly for the other side. Figure~\ref{figFalg} part (a) then checks that we have a Frobenius algebra with respect to this pairing and that the categorical adjoint of the product with respect to the pairing is the coproduct, and vice versa. Note that these categorical adjoints are the opposite of the usual vector space adjoints so as to avoid an unnecessary transposition, as for braided-Hopf algebras in \cite{Ma:alg}. (ii) Conversely, suppose that $A$ is a Frobenius algebra with bilinear $(\ ,\ )=\cup$ and inverse $g=\cap$. It is convenient to define $\Delta$ in terms of the bilinear and $\mu$ as shown in Figure~\ref{figFalg}, and $\eps=(\ ,1)$ as adjoint to 1. This implies that $\eps\mu=(\ ,\ )$, $\Delta(1)=g$ and $1=(\id\tens\eps)(g)$.  In the figure, we also recover $\mu$ in terms of $\Delta$ by the up-side-down formula and check the main identities for an F-algebra. Once these are established, we go back in the last line of part (b) and check that $\Delta$ is coassociative. From an identity in part (a) we also see that $\Delta,\mu$ are adjoint in the categorical sense. It is clear that this construction in inverse to the one in part (i), proving the equivalence. (iii) Part (c) of the figure makes it clear that the special property corresponds to $\mu(g)$ being a left and right identity for the product $\mu$, i.e. 1. (The same proof up-side-down means it is also equivalent to $(\ ,\ )\Delta=\eps$.)
 \endproof
 
 The significance of the F-algebra point of view is that all the identities for associativity, unity, coassociativity, counity, and the identities in  Definition~\ref{defFalg} take the form that all different ways to compose $\Delta,\mu, 1, \eps$ to give a connected planar graph with fixed $m$ legs in and $n$ legs out are equal. Here  $1$ is regarded as a node with no legs in and one leg out ($m=0,n=1$) and $\eps$ is regarded as a node with one leg in and no legs out ($m=1,n=0$). This leads iteratively to:
 
 \begin{corollary}\cite{MaRie} cf\cite{CD,Q} For a special  F-algebra, all compositions of $\Delta,\mu, 1, \eps$ corresponding to a connected planar graph with $m$ legs in and $n$ legs out can be obtained in the standard form of an iterated $m-1$-fold product and an iterated $n-1$-fold coproduct. By associativity and coassociativity, we can depict these as a single node with $m$-legs going in at the top and $n$ legs out at the bottom (a `spider'), in our case keeping the order of legs. 
 \end{corollary}

 Composing spiders along a contiguous subset of legs necessarily gives the spider with the remaining legs, which is also clear using coassociativity and associativity to isolate the branches to be contracted and using the special property iteratively. The empty morphism is the element 1 of the field, which must therefore equal $\eps(1)$ as the spider with $m=0,n=0$. (This is not really a connected graph but we formally include it.)  The same will apply up to an overall scale factor in the quasispecial case. The commutative and cocommutative case of primary interest in Computer Science  (and where we do not need to restrict to planar graphs, i.e., can allow crossings) is the celebrated spider theorem in \cite{CD}. The noncommutative case but with the Frobenius form symmetric occurs somewhat implicitly in \cite{LauPf}, and is also studied in \cite{Q}. A short proof in the fully general case is in \cite{MaRie}.

 \subsection{F-Hopf algebras}\label{secFHopf}

We recall that an algebra and coalgebra $H$ form a bialgebra if $\Delta,\eps$ are unital algebra homs (where $H\tens H$ has the tensor product algebra structure) and a Hopf algebra if there is additionally a map $S:H\to H$ such $\mu(S\tens\id)\Delta=1\eps=\mu(\id\tens S)\Delta$. See \cite{Ma} for more details. 
A single F-algebra can never be a bialgebra other than the trivial case where $A=k$ is the field (which we exclude). This is because $(h,h')=\eps(hh')$ for all $h,h'\in A$ and $g=\Delta(1)=1\tens 1$ need to be inverses to each other, which needs $1\eps$ to be the identity map on $A$. 

\begin{remark}\label{remunnorm}\rm In the Computer Science literature, a `bialgebra' often means what we will call an {\em unnormalised bialgebra} (or `scaled bialgebra'\cite{CD}), where we assume an algebra and a coalgebra and only that $\eps(1)$ is invertible and 
\[1\tens 1= \eps(1)\Delta(1),\quad \eps\tens\eps=\eps(1)\eps\mu,\quad \Delta\mu=\eps(1)\mu_{A\tens A}(\Delta\tens\Delta)\]
\[ \includegraphics{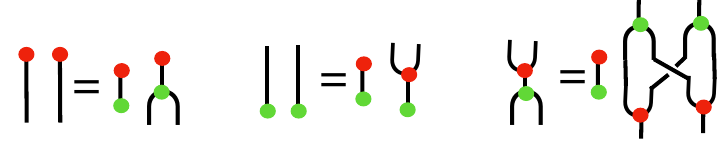}\]
There is no significance to the braid crossing at the moment (the operation so far is just the usual flip map). The unit and counit are represented by a univalent node as before. 
In this case, it is easy to see that $\lambda\Delta,\lambda^{-1}\eps$ form a usual bialgebra. An unnormalised Hopf algebra is an unnormalised bialgebra and an antipode $S$ obeying the same condition as for bialgebras. In this case $\lambda^{-2}S$ makes the associated bialgebra into a Hopf algebra in the usual sense. Hence it suffices to work in the usual normalised case. There is also no significance so far to colouring of the algebra and coalgebra nodes of a bialgebra other than to remind that they cannot belong to the same F-algebra.\end{remark}

\begin{definition}\cite{DD, CoD} An F-bialgebra is a pair of F-algebras on the same vector space such that for each F-algebra, the algebra of one and the coalgebra of the other forms a bialgebra. It is called an F-Hopf algebra if these bialgebras are Hopf algebras. It is called special if both F-algebras are special. 
\end{definition}

Here is where, because there are now two F-algebras, it is useful colour one red and the other green (in a black and white, printout, the green will appear lighter). Then the definition
says that they fit together e.g. as in Remark~\ref{remunnorm}, except that we stay in the usual (normalised) case.  Similarly with the colours swapped. If we are given one, say the red product/green coproduct Hopf algebra, as primary then we refer to the other as the associated Hopf algebra in the pair. 

\begin{figure}
 \[ \includegraphics[scale=1]{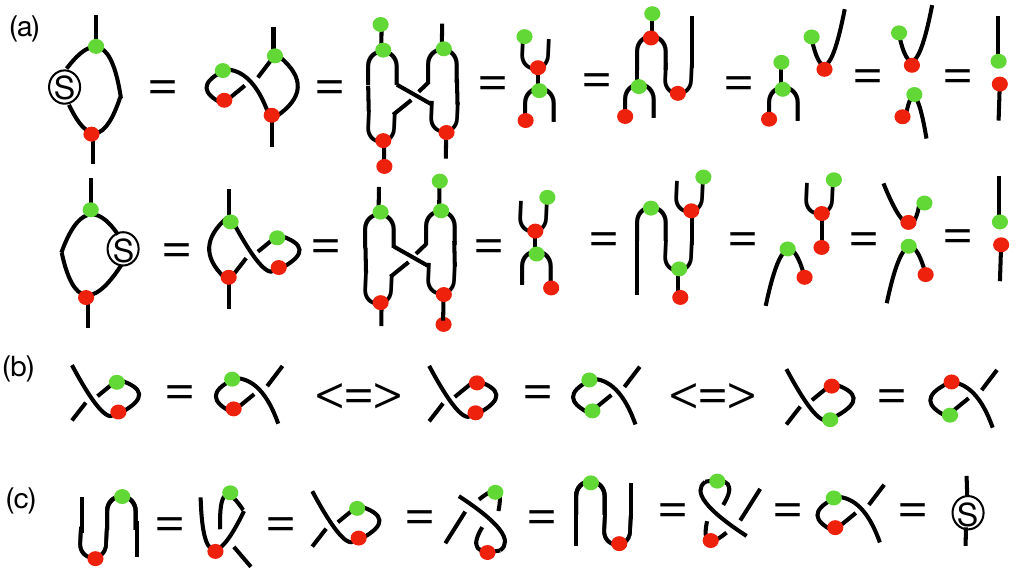}\]
 \caption{\label{figFant} (a) Proof of Proposition~\ref{propFant} that we obtain a Hopf algebra. (b) Three equivalent forms of the third condition. (c) Proof that the third condition holds in the symmetric case.}
 \end{figure}

\begin{proposition}cf\cite{DD}\label{propFant} If an $F$-bialgebra obeys  the three conditions 
\[ \includegraphics[scale=0.9]{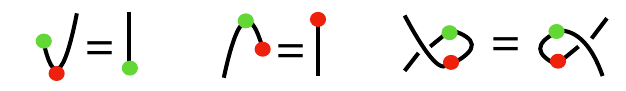}\]
then it is an $F$-Hopf algebra with antipode for the red product/green coproduct bialgebra given by the expressions in the third condition. The latter holds if the Frobenius forms are symmetric, in which case the antipode for the red product bialgebra and the green product bialgebra respectively can also be written as
\[ \includegraphics[scale=.8]{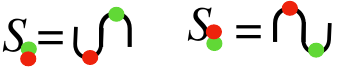}.\]
\end{proposition}
\proof The proof that we have an antipode for  the red product bialgebra is shown in Figure~\ref{figFant} (a) and is similar to  \cite{DD} but with a different $S$ as we do not assume here that our algebra is commutative. The left-right reflections of the first two stated conditions also hold by applying the unit or counit to the third stated condition. Next, applying the snake identities for duals to the third condition in the statement gives two equivalent versions in succession as shown in (b). One similarly has that the colour-reversal of the first two stated conditions hold. The colour-reversal of (a) is then the proof for the other Hopf algebra. 
Part (c) shows that the third stated condition holds if the red Frobenius form and the green metric are symmetric (with the latter equivalent to the green Frobenius form symmetric) and in this case it provides the second way to write the antipode for the red product bialgebra. Similarly for the  green product bialgebra . \endproof

The first two conditions here say that the categorical dual of the green unit with respect to the red Frobenius form is the green counit, and the green Frobenius dual of the red counit is the red unit.  This result is for motivation and in what follows we only take away the second special form of the antipode,  without the symmetry or other assumptions. In fact, we will assume $S$ to be the inverse of the one at the end of Proposition~\ref{propFant}, i.e. the stated form on the Hopf algebra with opposite coproduct.

In the following, we will start with one Hopf algebra $H$ and denote by $H^\medstar$ the categorical dual defined by nested adjoints (in the style of Figure~\ref{figFalg}(a), but now for the duality between $H$ and its dual space). Hence $H^\medstar$ is the op-algebra and op-coalgebra to the usual $H^*$, but the two are isomorphic as Hopf algebras via the antipode. We will assume that $H,H^\medstar$ as algebras are Frobenius. The latter means that $H$ is a {\em Frobenius coalgebra} in the sense of a coalgebra and an invertible metric $g$ such that
\[ (\Delta\tens\id)g=(\id\tens\Delta)g,\]
equivalent to  Definition~\ref{defFalg} by Lemma~\ref{lemFalg} read up side down.

\begin{figure}
 \[ \includegraphics[scale=1]{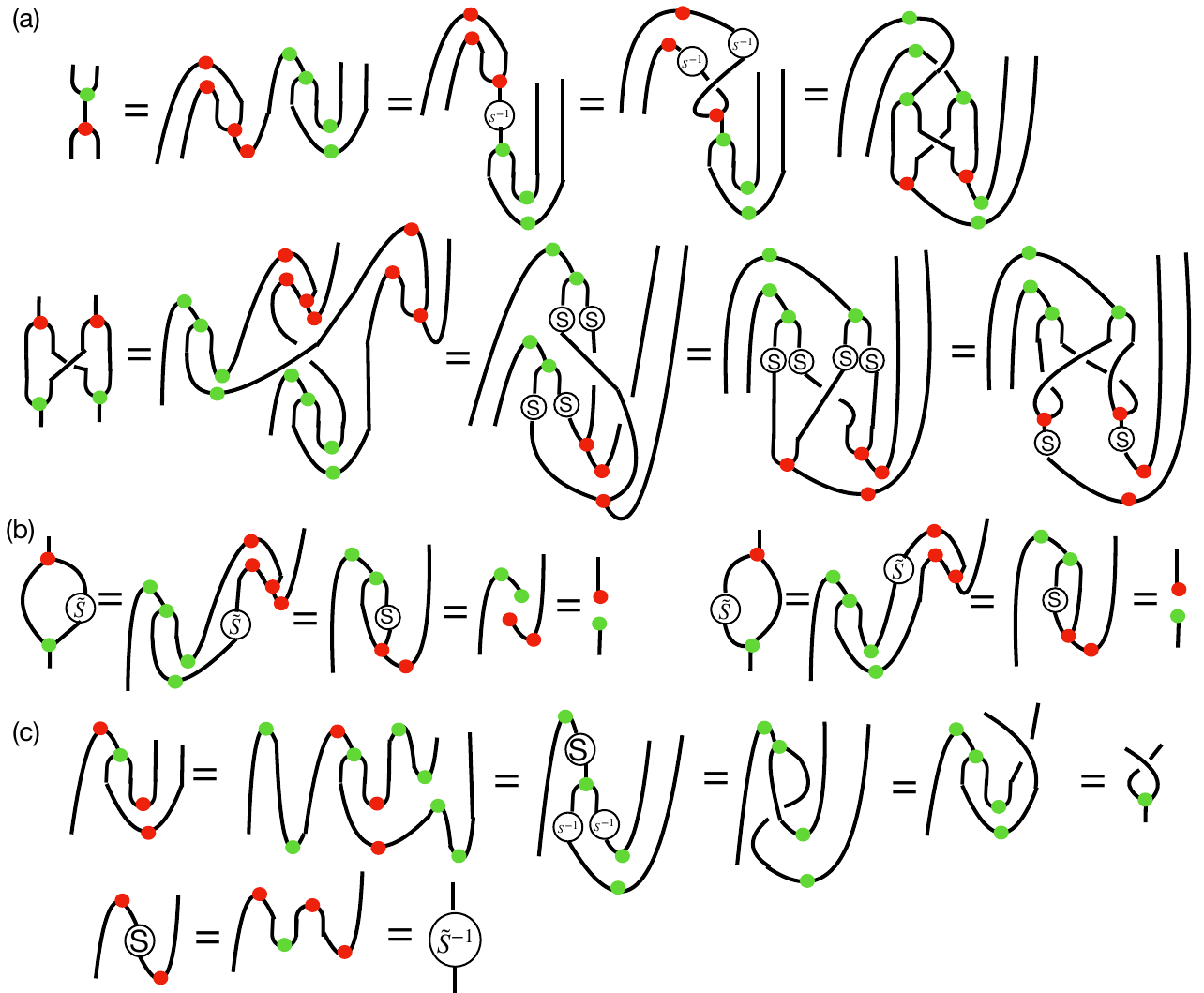}\]
 \caption{\label{figFhopf} Proof of Proposition~\ref{propFhopf}. Part (a) checks the bialgebra property and (b) the antipode properties.}
 \end{figure}

\begin{proposition}\label{propFhopf} Let $H=(\mu_{\color{red}\bullet},1_{\color{red}\bullet}, \Delta_{\color{green}\bullet},\eps_{\color{green}\bullet})$ be a Hopf algebra such that the algebra and coalgebra are Frobenius and the antipode $S$ is \includegraphics[scale=0.7]{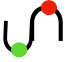}.
Then we have an F-Hopf algebra where the associated Hopf algebra $\tilde H=(\mu_{\color{green}\bullet},1_{\color{green}\bullet}, \Delta_{\color{red}\bullet},\eps_{\color{red}\bullet})$ has antipode  \includegraphics[scale=0.7]{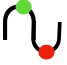}. Moreover, $\tilde H\isom H^{\medstar op}$ as Hopf algebras. 
\end{proposition}
\proof Given $H$ with red product and green coproduct we use the Frobenius structure on $H$ to define a dual red coproduct also on $H$. We use that the coalgebra is Frobenius to define a dual green product also on $H$. We then show in Figure~\ref{figFhopf} (a) that these fit together to form a bialgebra $\tilde H$, recognising $S^{-1}$ where $S$ denotes assumed antipode for the red product and green coproduct and using the (braided)-anti-algebra property of antipodes\cite{Ma} to take it through the product and cancel against the red metric, turning it green. In the second line {\em we start off with the wrong braid crossing}, which is of no significance at the moment, recognise $S$ and use the same (braided)-anti-algebra property to move them through the product. This then cancels against the red metric to turn it green and give the same expression as in the result of the first line after rearrangement. That the red counit (red-adjoint to the red unit) provides the counit of the bialgebra, etc are immediate and left to the reader. In part (b), we verify the proposed antipode $\tilde S$ for the green product/red coproduct Hopf algebra in terms of the antipode property for $S$. In (c) we check that $H^\medstar$ computed using (say) the red duality gives the product and antipode of $\tilde H^{op}$. 
\endproof
 
This gives a canonical F-Hopf algebra amplified under certain conditions from a single Hopf algebra. It provides a slightly different route to the following result in \cite{CD}, which also applies more generally for a Hopf algebra in a symmetric monoidal category such that integrals exist and give a nondegenerate bilinear form as for the vector space case on which we focus. 
 
 \begin{corollary}\label{corHH*}\cite{CD} Every finite-dimensional Hopf algebra $H$ and its dual $H^{\medstar op}$ form an $F$-Hopf algebra. This is quasispecial, with  $\mu_{\color{red}\bullet}(g_{\color{red}\bullet})\ne 0$ if and only if $H$ is semisimple and $\mu_{\color{green}\bullet}(g_{\color{green}\bullet})\ne0$ if and only if $H$ is cosemisimple (these being  equivalent in characteristic 0). 
 \end{corollary}
 \proof The algebra of any finite-dimensional Hopf algebra $H$ is necessarily Frobenius\cite{Par} with 
\[ (h,h')_{\color{red}\bullet}=\int hh',\quad g_{\color{red}\bullet}=(\id\tens S) \Delta\Lambda,\quad\mu_{\color{red}\bullet}(g_{\color{red}\bullet})=\eps(\Lambda)\]
\[ (\int\tens \id)\Delta= 1\int,\quad  \int\Lambda=1,\quad h\Lambda =\Lambda\eps(h)\]
for all $h,h'\in H$. The (red) product of $H$ is omitted and $\Delta$ is its (green) coproduct. The map $\int: H\to k$ here is  a `right-invariant integral' on $H$ and is unique up to scale. The element $\Lambda$ is a left integral in $H$ and corresponds to a left-invariant integral on $H^*$. Both are unique up to scale and can be normalised as stated. We use the compact notation $\Delta h=h_1\tens h_2$ (with more numbers for iterated coproducts and sum of such terms understood; alternatively this can be done with strings, given later). We check the metric inversion identities,
\begin{align*} (\int h\Lambda_1)S\Lambda_2&=(\int h_1\Lambda_1) (S(h_2 \Lambda_2))h_3=(\int h_1\Lambda)S(1)h_2=\eps(h_1)h_2\int\Lambda=h\\
\Lambda_1\int( (S\Lambda_2) Sh)&=(Sh_1)h_2\Lambda_1\int S(h_3\Lambda_2)=(Sh)\Lambda_1 \int S\Lambda_2\\
&=(Sh)S^{-1}(S\Lambda)_2\int (S\Lambda)_1=(Sh)S^{-1}(1)\int S\Lambda=Sh \int S\Lambda\end{align*}
for all $h\in H$. Here, computing the following middle expression two ways by the results just found,
\[ \Lambda_1'\tens S\Lambda_2'=\Lambda'_1\int (S\Lambda'_2)\Lambda_1\tens S\Lambda_2=\Lambda_1 \int S\Lambda\tens S\Lambda_2\]
tells us that $\int S\Lambda=1$, as required. Here $\Lambda'$ is another copy of $\Lambda$.  These calculations can also be done with diagrams, which we defer to Section~\ref{secbraZX}. This gives the red Frobenius form, i.e. on $H$ as an algebra. We have $\mu_{\color{red}\bullet}(g_{\color{red}\bullet})=\Lambda_1S\Lambda_2=1\eps(\Lambda)$, so this is quasispecial. A property of integrals is that  $\eps(\Lambda)\ne 0$ if and only if the algebra is semisimple. 

Moreover, the same data $\int,\Lambda$ can be viewed as $\int\in H^\medstar$ and $\Lambda:H^\medstar\to k$ and  make this similarly Frobenius, i.e. the coalgebra of $H$  Frobenius. We have explicitly
\[ (h,h')_{\color{green}\bullet}=\int (Sh)h'=(Sh,h')_{\color{red}\bullet},\quad g_{\color{green}\bullet}=\Delta\Lambda=(\id\tens S^{-1})g_{\color{red}\bullet}\]
from which it is clear that these are also inverse to each other, since the previous pair were, and that $S$ has the form needed in Proposition~\ref{propFhopf}. We have $\mu_{\color{green}\bullet}(g_{\color{green}\bullet})=\Lambda\int 1$, so as a Frobenius algebra, this too is quasispecial with nonzero value if and only if $H^\medstar$ is semisimple. The latter in characteristic 0,  is equivalent to $H$ semisimple and to $S^2=\id$ due to results of Larson and Radford, see \cite[Thm~3.14]{Sch} for an exposition. 
\endproof

For an example, we recall\cite{Ma} that if $X$ is a group (not necessarily Abelian) then its associated group Hopf algebra $kX$  has basis given by elements of $X$ with the product extended linearly and $\Delta x=x\tens x, \eps(x)=1, Sx=x^{-1}$ for all $x\in X$. Elements, as here, where the coproduct is diagonal are called {\em group-like}\cite{Ma} (or set-like in the computer science literature). If $X$ is finite then the associated function Hopf algebra  $k(X)$ has basis $\{\delta_x\}$ of $\delta$-functions with $\delta_x\delta_y=\delta_{x,y}\delta_x$, $\Delta\delta_x=\sum_{yz=x}\delta_y\tens\delta_z$, $\eps\delta_x=\delta_{x,e}$ and $S\delta_x=\delta_{x^{-1}}$ for all $x\in X$. Here $e\in X$ is the group identity.

\begin{example}\label{exCX}\rm  If $X$ is a finite group then $H=kX$ with product $x{\color{red}\bullet}y=xy$ and coproduct $\Delta_{\color{green}\bullet}x=x\tens x$ amplifies to an $F$-Hopf algebra by Proposition~\ref{corHH*}. The integral is $\int x=\delta_{x,e}$ and the integral element is $\Lambda=\sum_{x\in X}x$. The Frobenius form on $H$ as an algebra, its inverse, and the associated coproduct are
\[ (x,y)_{\color{red}\bullet}=\delta_{x^{-1},y},\quad g_{\color{red}\bullet}=\sum_{x\in X} x\tens x^{-1},\quad \Delta_{\color{red}\bullet}x=\sum_{y,z\in X\, |\, yz=x}y\tens z,\quad\eps_{\color{red}\bullet}x=\delta_{e,x}\]
and $\mu_{\color{red}\bullet}(g_{\color{red}\bullet})=|X|$, which is quasi-special. The Frobenius form on $H$ as a coalgebra, its inverse and the associated product are
\[(x,y)_{\color{green}\bullet}=\delta_{x,y},\quad g_{\color{green}\bullet}=\sum_{x\in X} x\tens x,\quad x_{\color{green}\bullet}y= \delta_{x,y}x,\quad 1_{\color{green}\bullet}=\sum_{x\in X}x \]
from which we see that the green algebra/red coalgebra has the structure of $H^*$ on matching $x$ to $\delta_x$ there. We also see that the green F-algebra is special.
 \end{example}

 The work \cite{CoD} considers the simplest noncommutative noncocommutative example, the Taft algebra. Going beyond this, we will now look at the reduced Drinfeld-Jimbo quantum group $u_q(su_2)$ at $q$ a root of unity. There are different conventions for this and we use essentially the ones in the recent work\cite{AziMa} which is noncommutative and noncocommutative even when $q=-1$, when understood correctly. 
 
 \begin{example}\label{exsl2}\rm Let $q$ be a primitive $n$-th root of unity and let $H=u_q(sl_2)$ as generated by $E,F,K$ and 
\[ E^n=F^n=0,\quad K^n=1,\quad KEK^{-1}=q^{-1}E,\quad KFK^{-1}=qF,\quad [E,F]={K-K^{-1}\over q-q^{-1}}\]
\[ \Delta K=K\tens K,\quad \Delta E=1\tens E+E\tens K,\quad \Delta F=F\tens 1+ K^{-1}\tens F,\quad\eps K=1,\quad \eps E=\eps F=0\]
with antipode (and the quasitriangular structure, needed later),
\[ SK=K^{-1},\quad SE=-EK^{-1},\quad SF=-KF,\quad \CR={1\over n}\sum_{r,a,b=0}^{n-1}{(-1)^r(q-q^{-1})^rq^{-ab}\over [r]_{q^{-1}}!}F^rK^a\tens E^r K^b\]
where $[m]_q:=(1-q^m)/(1-q)$ denotes a `q-integer'. When $q=-1$ we omit the $q-q^{-1}$ factors in the relations and in $\CR$. The right integral on $H$ and left integral in $H$ are
\[ \int K^iF^j E^k=\delta_{1,i}\delta_{n-1,j}\delta_{n-1,k},\quad \Lambda=\Lambda_KF^{n-1}E^{n-1};\quad \Lambda_K=\sum_{r=0}^{n-1}K^r\]
so that $\int\Lambda=1$. One also has that $S\Lambda=\Lambda$. We compute the associated $\tilde H$ product and coproduct
\[ h{\color{green}\bullet}h'=h_1\int ((Sh_2)h'),\quad  \Delta_{\color{red}\bullet}h=\Lambda_1\tens (S\Lambda_2)h\]
for all $h,h'\in u_q(sl_2)$ can then be computed and must be isomorphic to $c_q[SL_2]^{op}$ for the appropriate reduced quantum group  dual to $u_q(sl_2)$. We have $\mu_{\color{red}\bullet}(g_{\color{red}\bullet})=\mu_{\color{green}\bullet}(g_{\color{green}\bullet})=0$. 

We compute these products for the $n=2$ case where $q=-1$ and\cite{AziMa}
\[K^2=1,\quad E^2=F^2=0,\quad  [E,F]=0,\quad  \{E,K\}=\{F,K\}=0,\]
\[\Delta K=K\tens K, \quad \Delta E=E\tens K+1\tens E,\quad \Delta F=F\tens 1+K\tens F\]
\[  \eps K=1,\quad \eps E=\eps F=0,\quad SK=K,\quad  SE=KE,\quad  SF=-KF\]
\[\Lambda=(1+K)EF,\quad \CR=(1\tens 1-F\tens E)\CR_K,\quad \CR_K={1\over 2}(1\tens 1+K\tens 1+1\tens K-K\tens K).\]
This is the initial Hopf algebra with red algebra and green coalgebra understood. 
From these we work out the values of $(\id\tens S)\Delta$ on our monomial basis. In particular,
\[ g_{\color{red}\bullet}=\Lambda_1\tens S\Lambda_2=F\tens KE-KE\tens F+ EF\tens K+K\tens EF+E\tens KF-KF\tens E+KEF\tens 1+1\tens KEF\]
We let 
\[ x=(1+K)E,\quad y=(1+K)F,\quad t=(K-1)EF\]
and find for the associated Hopf algebra
\[ \Delta_{\color{red}\bullet} t=t\tens t,\quad \Delta_{\color{red}\bullet} x=x\tens t+\Lambda\tens x,\quad \Delta_{\color{red}\bullet} y=y\tens t+\Lambda\tens y\]
\[ x{\color{green}\bullet}x=y{\color{green}\bullet}y=0,\quad t{\color{green}\bullet}t=\Lambda\]
\[ x{\color{green}\bullet}y=-y{\color{green}\bullet}x=1+K,\quad x{\color{green}\bullet}t=-t{\color{green}\bullet}x=(K-1)E,\quad y{\color{green}\bullet}t=-t{\color{green}\bullet}y=(1-K)F.\]
Note that $1_{\color{green}\bullet}=\Lambda$ is the unit element. In fact, $u_q(sl_2)$ at $q=-1$ is both self-dual and anti-self dual, and indeed the associated Hopf algebra here is isomorphic to $u_q(sl_2)$  if we identify
\[ 1\mapsto 1_{\color{green}\bullet},\quad K\mapsto t,\quad E\mapsto x,\quad F\mapsto y{\color{green}\bullet}t.\]
\end{example}

\section{ $*$-algebra and Hadamard forms} \label{secstar}

Next, we recall the usual notion of a Hopf $*$-algebra. This means a Hopf algebra over $\C$ where the algebra is a $*$-algebra (i.e., it has an antilinear anti-algebra involution $*$) and $\Delta,\eps$ commute with $*$ (which on $\C$ is complex conjugation). In this case, one has that $(S*)^2=\id$. For a group algebra $\C X$, we have $x^*=x^{-1}$ and for a function Hopf algebra $\C(X)$ we have $f^*(x)=\overline{f(x)}$ where the bar denotes complex conjugation of the value. This is the standard notion but will not be the most relevant. Instead, there is another which we called a {\rm flip-Hopf} $*$-algebra in \cite{Ma, BegMa}, where
\[ \Delta * =\dagger\Delta,\quad \overline{\ }\circ\eps=\eps\circ *,\quad S*=*S\]
with the `hermitian conjugation' 
\[ \dagger:={\rm flip}(*\tens *)\]
on any tensor product space over $\C$. The categorical picture behind these ideas is that of a {\rm bar} category\cite{BegMa:bar,BegMa} which will be particularly useful in the braided case.  For any flip-Hopf $*$-algebra, we define the associated antilinear Hopf algebra automorphism
\[ \theta:= S*,\]
which squares to $S^2$. 

\subsection{Unimodular setting}\label{secuni}

We start with the simplest setting, which is adequate for finite group examples and simpler quantum groups.

\begin{definition}\label{starform} A Frobenius $*$-algebra is a Frobenius algebra $A$ over $\C$ which is a $*$-algebra and for which $\overline{(a,b)}=(b^*,a^*)$, or equivalently $g^\dagger=g$ for the metric. 
\end{definition}
The metric property here is similar to the notion of a real metric in noncommutative geometry \cite{BegMa} and the proof that this is equivalent to the property of (\ ,\ ) is parallel to the proof there. In this case, it makes sense to define a  sesquilinear bilinear form on $H$ by
\[ \<h|h'\>=(h^*,h')\]
which, under our assumption makes $A$ into a Hilbert space (upon completion in the infinite-dimensional case). In the group algebra Example~\ref{exCX} we have $\<x|y\>=\delta_{x,y}$ for the (red) inner product there, i.e. the group basis elements become an orthonormal basis.

\begin{lemma}\label{starfrob} A $*$-algebra $A$ is Frobenius if and only if the corresponding F-algebra obeys \rm
\[ \Delta * =\dagger \Delta,\quad \overline{\ }\circ\eps=\eps\circ *.\]
\end{lemma}
\proof Here $\Delta a=ga$ and we use Lemma~\ref{gcentral} to deduce that $\Delta (a^*)= g a^*=a^*g=\dagger(g^\dagger a)=\dagger(ga) =\dagger\Delta a$. We used that the metric obeys $g^\dagger=g$. We also have $\eps(a^*)=(1,a^*)=\overline{(a,1)}=\overline{\eps(a.1)}=\overline{\eps(a)}$. \endproof

We equally define a Frobenius $*$-coalgebra as a Frobenius coalgebra over $\C$ where the coalgebra obeys the conditions in Lemma~\ref{starfrob} and the Frobenius structure obeys the conditions in Definition~\ref{starform}, in which case the associated algebra is a $*$-algebra. In either case, we say that the associated F-algebra is an F-$*$-algebra.  

\begin{corollary}\label{corbasicHH*}(Unimodular case.) If $H$ is a flip-Hopf $*$-algebra and $\int *=\overline{\ }\circ\int, \Lambda^*=\Lambda$ in Corollary~\ref{corHH*}  then both associated F-algebras are F-$*$-algebras and the associated Hopf algebra  in the induced F-Hopf algebra is also a flip-Hopf $*$-algebra. 
\end{corollary}
\proof Here $(h'{}^*,h^*)_{\color{red}\bullet}=\int(hh')^*=\overline{\int(h,h')_{\color{red}\bullet}}$ in Corollary~\ref{corHH*} under our assumption on $\int$. Similarly, $\dagger g_{\color{green}\bullet}=\dagger\Delta\Lambda=\Delta(\Lambda^*)=g_{\color{green}\bullet}$ for the Frobenius coalgebra in Corollary~\ref{corHH*}. Hence the red algebra and the green coalgebra are both Frobenius as *-(co)algebras. The associated red coproduct therefore obeys the condition in Lemma~\ref{starfrob} and the green product is a *-algebra. Hence the Hopf algebra that these form by Corollary~\ref{corHH*} is a flip-Hopf $*$-algebra. \endproof

We say in this case that we have an example of an F-Hopf $*$-algebra. 

\begin{example}\rm \label{exCXstar} Let $H=\C X$ for a finite group $X$ with $x^*=x^{-1}$. This is a Hopf $*$-algebra but we regard it as a flip-Hopf $*$-algebra since it is cocommutative. We have $\int x^*=\delta_{e,x^{-1}}=\delta_{e,x}=\int x=\overline\int x$ and $\Lambda^*=\sum_{x\in X}x^{-1}=\Lambda$ as required. The induced red coproduct in Example~\ref{exCX}  obeys $\Delta_{\color{red}\bullet}(x^*)=\sum_{yz=x^{-1}}y\tens z=\dagger\sum_{yz=x^{-1}}z^{-1}\tens y^{-1}=\dagger\Delta_{\color{red}\bullet} x$ as it had to by Corollary~\ref{corbasicHH*}. Likewise the induced green product obeys $y^*{\color{green}\bullet}x^*=y^{-1}{\color{green}\bullet}x^{-1}=\delta_{x,y}x^{-1}=(x{\color{green}\bullet}y)^*$, so this is a $*$-algebra as it had to be for an F-Hopf $*$-algebra.\end{example}

Note that $\C(X)$ with its usual $*$-algebra structure is not a flip-Hopf *-algebra when $X$ is nonAbelian. For that, we should take $\delta_x^*=\delta_{x^{-1}}$. This then matches up with the red coproduct and green product just found. 

\begin{example}\label{exflipsl2}\rm Let $H=u_q(sl_2)$ at $q=-1$ as in Example~\ref{exsl2}. This forms a flip-Hopf $*$-algebra with 
\[ K^*=K,\quad E^*=F,\quad F^*=E;\quad \theta(K)=K,\quad \theta(E)=-KF,\quad \theta(F)=KE.\]
We also have $\int S= \int$, $S\Lambda=\Lambda$ (the Hopf algebra is unimodular in both aspects) and $\int h^*=\overline{\int h}$,\ $\Lambda^*=\Lambda$ as required. Using the same underlying $*$-operation, we have
\[ t^*=FE(K-1)=t,\quad x^*=F(1+K)=(1-K)F=y{\color{green}\bullet}t,\quad y^*=E(1+K)=(1-K)E=t\circ x\]
which indeed makes the green product a $*$-algebra also. Similarly, for example
\begin{align*}\Delta_{\color{red}\bullet}x^*&=\Delta_{\color{red}\bullet}(y{\color{green}\bullet}t)=(\Delta_{\color{red}\bullet}y){\color{green}\bullet}(\Delta_{\color{red}\bullet}t)=(y\tens t+\Lambda\tens y){\color{green}\bullet}(t\tens t)=y{\color{green}\bullet}t\tens\Lambda+t\tens y{\color{green}\bullet}t\\
&=x^*\tens\Lambda+t\tens x^*=\dagger\Delta_{\color{red}\bullet} x \end{align*}
as it must for a flip-Hopf $*$-algebra. In fact, this is isomorphic as a flip-Hopf $*$-algebra to $u_q(sl_2)$ by the map in Example~\ref{exsl2}. 
\end{example}

On the other hand, the assumptions we made  are not typical because in general $\int*$ is an (antilinear) left invariant integral and can not be equated to $\overline{\ }\circ\int$, and $\Lambda^*$ is a right integral element and can not be equated to $\Lambda$. This is because a typical quantum group will not be unimodular in the sense that the left and right integrals will not be proportional. Note that the Frobenius form, when it exists, is not unique. Hence it is also possible that we could stay in the setting of Definition~\ref{starform} and Lemma~\ref{starfrob}, which seems reasonable from a noncommutative geometry point of view, but with the Frobenius form and its inverse not given by integrals. However, we do not have any general results about this more general possibility. One could also explore a version of the theory adapted to ordinary Hopf $*$-algebras, which is more relevant to usual quantum geometry and applicable to $U_q(su_2)$ and $\C_q[SU_2]$ but not to the reduced finite-dimensional versions at $q$ a root of unity. 

\subsection{General Hopf algebra case based on integrals} \label{secgenstar}

The less special assumptions we will take for the properties of the integrals on a flip-Hopf $*$-algebra under $*$ are
\begin{equation}\label{intstar} \int h^*= \overline{\int S h},\quad \Lambda^*=S\Lambda\end{equation}
These have to hold up to scale as each space of integrals is 1-dimensional, so it is not unreasonable to assume them as part of the normalisations. This does, however, imply that $\int\circ S^2=\int$ and $S^2\Lambda=\Lambda$

\begin{proposition}\label{propflipHH*} Let $H$ be a flip-Hopf $*$-algebra with  integrals in Corollary~\ref{corHH*} obeying the condition (\ref{intstar}) and suppose that the antipodes on $H$ and on associated Hopf algebra of the induced  F-Hopf algebra square to the same. Then the associated Hopf algebra is a flip-Hopf $*$-algebra with the same associated antilinear automorphism. \end{proposition}
\proof The issue here is that $\int,\Lambda$ are not assumed to obey the stricter reality properties in Corollary~\ref{corbasicHH*}, with the result that the associated green product is not a $*$-algebra with the same $*$ as $H$, and likewise the red coproduct does not skew-commute with $*$. We will resolve this by introducing a new green $*$-operation $\color{green}*$ making this into a flip-Hopf $*$-algebra. First note that the first of (\ref{intstar}) says that $\overline{\int h}=\int\theta(h)$. Hence
\begin{align*}\theta(h {\color{green}\bullet}h')& = \theta(h_1{})\overline{(Sh_2)h'}= \theta(h_1{})\int \theta(Sh_2)\theta(h')=\theta(h{})_1\int (S\theta(h)_2)\theta(h')=\theta(h){\color{green}\bullet}\theta(h').\end{align*}
Next, the second of (\ref{intstar}) implies that 
\begin{align*} g^\dagger&=(S\Lambda_2)^*\tens\Lambda_1{}*=S(\Lambda_2{}^*)\tens\Lambda_1{}*=S\Lambda^*{}_1\tens\Lambda^*{}_2\\
&=S(S\Lambda)_1\tens (S\Lambda)_2=S^2\Lambda_2\tens S\Lambda_1=(S\tens S){\rm flip}(g)\end{align*}
which is equivalent to $(\theta\tens\theta)(g)=g$. Consequently,
\[ (\theta\tens\theta)\Delta_{\color{red}\bullet}h=(\theta\tens\theta)(gh)=g\theta(h)=\Delta_{\color{red}\bullet}\theta(h)\]
This $\theta$ is {\em also} a bialgebra and hence Hopf algebra automorphism on the associated green product/red coproduct Hopf algebra. If we let $\tilde S$ denote the antipode for the associated green algebra/red coalgebra $\tilde H$ then $\tilde S\theta=\theta \tilde S$. We now see that 
${\color{green}*}:=\tilde S^{-1}\theta$ is an antilinear anti-algebra and anti-coalgebra map and commutes with $\tilde S$. Hence 
\[ {\color{green}*}^2=\tilde S^{-1}\theta \tilde S^{-1}\theta=\theta^2 S^{-2}_{\color{green}\bullet}=S^2 S^{-2}_{\color{green}\bullet}=\id\]
if and only if $S^2=S^2_{\color{green}\bullet}$. Under this additional hypothesis, we see that ${\color{green}*}$ makes the associated Hopf algebra into a flip-Hopf $*$-algebra. 
 \endproof

In this case, we say that we have a general F-Hopf $*$-algebra, where we have a green $*$-operation for the green product, not necessarily the same as the initial one for the initial (red) product. The notion of a Frobenius $*$-algebra also gets modified. We saw that now the natural property is that $(\theta\tens\theta)(g_{\color{red}\bullet})=g_{\color{red}\bullet}$ and one similarly has 
\[\overline{(h,h')_{\color{red}\bullet}}=\overline{\int hh'}=\int\theta(h)\theta(h')=(\theta(h),\theta(h'))_{\color{red}\bullet}\]
under our assumptions. The same on the green side, noting that our definition of ${\color{green}*}$ is such that the green $\theta$ automorphism is the same as the initial one on $H$. This therefore departs from the more natural setting from the Frobenius $*$-algebra point of view in the preceding section. 

\begin{example}\rm  We make $H=u_q(sl_2)$ in Example~\ref{exsl2} for $q$ a primitive $n$-th root of unity into a  flip-Hopf $*$-algebra\cite{Ma} with 
\[ K^*=K^{-1},\quad E^*=F,\quad F^*=E;\quad \CR^\dagger=\CR^{-1}.\]
The right integral on $H$ and left integral in $H$ obey the conditions (\ref{intstar}) (in fact, one has that  $\Lambda^*=\Lambda=S\Lambda$ is unimodular, but not $\int$). These are straightforward calculations which we omit, focusing on the additional condition in Proposition~\ref{propflipHH*}. In lieu of a formal proof, we do a detailed calculation to illustrate the method. Thus, from Corollary~\ref{corHH*}, 
\[ \tilde Sh=\Lambda_1\int \Lambda_2 h=\sum_j(\id\tens\int)\Delta (K^j F^{n-1}E^{n-1})h\]
which on a monomial element $h$ picks out one term of the coproduct such that the powers of $F,E$ in the right hand factor  add up with the powers in $h$ to $n-1$ (else the integral will kill the term). For example $h=K^a F$ means the only relevant terms are 
\[ \Delta E^{n-1}= (E\tens K+1\tens E)^{n-1}=1\tens E^{n-1}+...\]
\[ \Delta F^{n-1}=(F\tens 1+K^{-1}\tens F)^{n-1}=...+[n-1]_{q^{-1}} FK^{-(n-2)}\tens F^{n-2}+...\]
on expanding in powers of $E,F$ in the right factors and using the commutation relations of $u_q(sl_2)$ to collect terms (an instance of the $q$-binomial formula\cite{Ma}).  Hence,
\begin{align*} \tilde S(K^aF)&=\sum_j(\id\tens\int)\left( (K^j\tens K^j)([n-1]_{q^{-1}} FK^{-(n-2)}\tens F^{n-2})(1\tens E^{n-1} K^a F)\right)\\
&=\sum_j [n-1]_{q^{-1}}K^j F K^{2}\int K^j F^{n-2} E^{n-1}K^aF\\
&=\sum_j [n-1]_{q^{-1}} q^{a-2} K^{2+j}F\int  K^{a+j} F^{n-1} E^{n-1}\\ &=\sum_j [n-1]_{q^{-1}} q^{a-2} K^{2+j}F\delta_{1,a+j}=[n-1]_{q^{-1}} q^{a-2} K^{3-a}F=-q^{a-1}K^{3-a}F.
\end{align*}
We commuted $E,F$ inside the integral as their commutator has lower degree in $E,F$. Hence $S^2_{\color{green}\bullet}(K^aF)=q^{a-1+ (3-a)-1}K^{3-(3-a)}F=q K^aF$. Meanwhile, $S(K^aF)=-KFK^{-a}=-q^a K^{1-a}F$ so that $S^2(K^aF)=q^{a+1-a}K^{1-(1-a)}= qK^aF$ also, so $S^2_{\color{green}\bullet}=S^2$ on these elements. The proof for the general case is rather more complicated, but follows the same method using the $q$-binomial formula and extracting the relevant terms for the coproduct. Thus, the conditions of Proposition~\ref{propflipHH*} apply and we have an F-Hopf $*$-algebra. The associated Hopf algebra in this example is presumably isomorphic to a reduced version of the flip-Hopf $*$-algebra  $\C_q[SU_2^{fl}]$ in \cite[Ex.~2.113]{BegMa}.  \end{example}

\subsection{Phases} Another important ingredient for ZX calculus is the notion of phases. These are defined in \cite{CD,DD} as map $A\to A$ on a Frobenius algebra $A$ such that applied to either leg of the product gives the same result as applying the map after the product. It is shown in the commutative setting that these take the form of left multiplication by elements of $A$. This map can also be applied to coproducts leading to the decoration of any spider by this operation applied to any leg.

In our setting with $A$ a possibly noncommutative Frobenius $*$-algebra, we let $Z(A)$ be the centre of $A$  and  
\[ P(A)=\{\alpha\in Z(A)\ |\  \alpha^*\alpha=\alpha\alpha^*=1\}\]
i.e. of invertible elements such that $\alpha^*=\alpha^{-1}$. This forms an Abelian group (the `group of phases'). The corresponding operation is $L_\alpha(a)=\alpha a$ as before, which clearly can be applied to any leg of a composite product. 

\begin{lemma} If $A$ if an F-$*$-algebra then the action of $\alpha\in P(A)$ can be taken through the coproduct to either leg.
\end{lemma}
\proof This is immediate using the relevant form of the coproduct. For example, $(\id\tens L_\alpha)\Delta a=g^1\tens \alpha g^2 a=g^1\tens g^2\alpha a=\Delta(L_\alpha a)$ and $(L_\alpha\tens \id)\Delta a=\alpha g^1\tens  g^2 a=g^1\tens g^2\alpha a=\Delta(L_\alpha a)$ using centrality of $g$ (or using the other form $\Delta a=ag$.) \endproof

It follows, as in \cite{CD,DD}, that we get the same result applying $L_\alpha$ to any leg of a spider, and hence that we can denote such an $\alpha$ in a neutral way inside the node of the spider. When combining spiders, it is also clear that phases multiply using the group structure of $P(A)$ inherited from the product of $A$. In \cite{CD,DD}, there is also a discussion of group-like (`set-like') elements and classical values which we omit. 

\subsection{Hadamard gates and self-duality}

Another ingredient of $ZX$ calculus is the Hadamard gate, which for the standard 1 qubit realisation using $\C\Z_2$ in the introduction, is clearly just the Fourier transform. This observation was used in \cite{GogZen} to propose Fourier theory on internal groups in their categorical setting of ZX calculus. Here we give a different take on the role of this gate in our algebraic setting.

 Note first that Fourier transform works canonically on any finite-dimensional Hopf algebra\cite{Ma,BegMa}. If $H$ is a Hopf algebra with right-invariant integral as above, then the Hopf algebra Fourier transform is an invertible linear map
\begin{equation}\label{fou} \CF: H\to H^\medstar,\quad \CF(h)=\sum_a f^a(\int e_a h),\end{equation}
where $\{e_a\}$ is  a basis of $H$ and $\{f_a\}$ a dual basis. We have seen that in the F-Hopf algebra case, the associated (green product, red coproduct) Hopf algebra is isomorphic to $H^\medstar$ as a vector space and indeed if we use the metric $g_{\color{red}\bullet}$ to refer $\CF$ back to a map $H\to H$ then it is just the identity map. Thus, the canonical Fourier transform itself is not the content of the Hadamard gate. 

In fact, Fourier transform as encountered on $\R^n$ or on $\Z_n$ in terms of the group algebra is an operator $H\to H$ after the canonical Fourier transform is combined with the fact that these groups are self-dual. 
We propose this self-duality of some Hopf algebras as an approach to a Hadamard gate, but different from another point of view that it should map the green spiders to red spiders, which seems important for ZX calculus. We start with the latter. 

\begin{definition}\label{defhad}  A Type 1 {\em Hadamard form} on an $F$-bialgebra $H$ is a non-degenerate bilinear form $\Theta$ on $H$ such that (a)-(b) of
\[ \includegraphics{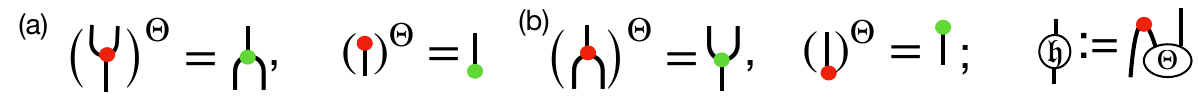}\]
hold for adjunction with respect to $\Theta$.  $\Theta$ is called {\em quasi-Hadamard} if these hold up to constant scale factors. \end{definition}

We also show the definition of a corresponding `Hadamard gate' $\ch:H\to H$, which we will use throughout the section (also for the other types below). In the case where $(\ ,\ )_{\color{red}\bullet}=\int\circ\mu_{\color{red}\bullet}$, 
 \begin{equation}\label{hadfou} \ch^{-1}(h)=\Theta^1(\Theta^2,h )_{\color{red}\bullet}=\Theta^1\int \Theta^2 h\end{equation}
 is the Hopf algebra Fourier transform (\ref{fou}) with the associated `metric' $\Theta=\Theta^1\tens\Theta^2\in H\tens H$ in the role of `exponential'  used to map $H^\medstar\to H$ by evaluation against its second leg (i.e., used to provide the required duality).

\begin{figure}
 \[ \includegraphics[scale=1.1]{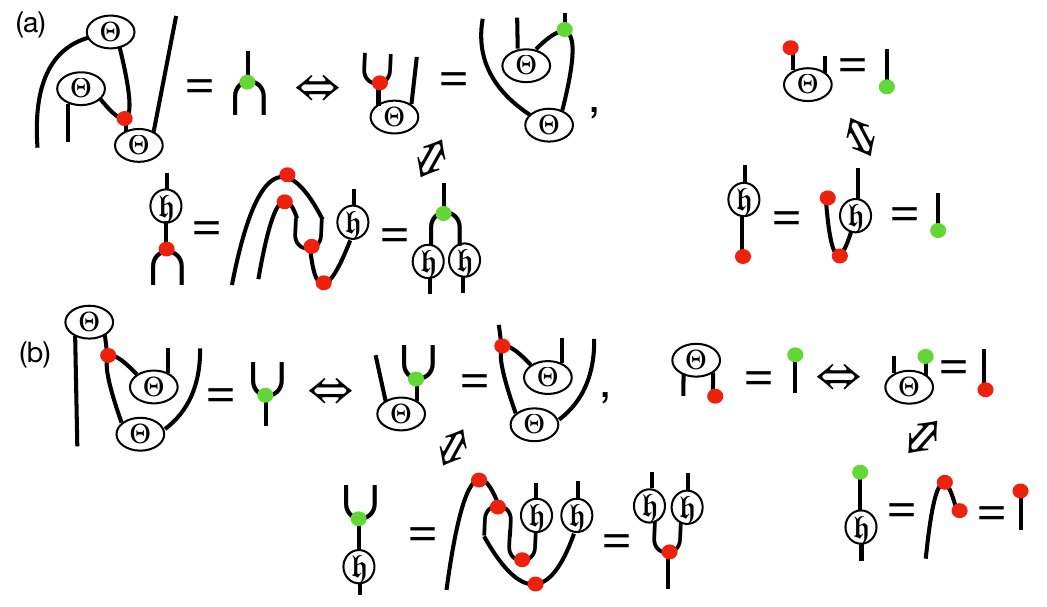}\]
 \caption{\label{fighadpf} Definition (a) and (b) of a Hadamard form on a Hopf algebra and equivalent properties of $\ch$.}
 \end{figure}

\begin{lemma} Let $H$ be an F-bialgebra or F-Hopf algebra and $\Theta$ a bilinear form on $H$. The following are equivalent.

(1) $\Theta$ is a Type 1 Hadamard form.

(2) $\ch$ is an isomorphism from the green F-algebra to the red F-algebra,
\[\mu_{\color{green}\bullet}=\ch^{-1}\circ\mu_{\color{red}\bullet}\circ(\ch\tens\ch),\quad \Delta_{\color{green}\bullet}=(\ch^{-1}\tens\ch^{-1})\circ\Delta_{\color{red}\bullet}\circ\ch,\quad  1_{\color{green}\bullet}=\ch^{-1}(1_{\color{red}\bullet}),\quad \eps_{\color{green}\bullet}=\eps_{\color{red}\bullet}\circ\ch\]
\end{lemma}
\proof The content of the Hadamard form assumption is shown explicitly in Figure~\ref{fighadpf}.   The proof is then straightforward but necessary to be sure that all relevant conventions match up, as shown in parts (a)-(b) of the figure. In (a) we cover the adjunction of the red algebra to obtain the green coalgebra, and in (b) the same for the algebra. In each case, we start with (1) the adjunction as a `rotation' using the bilinear form $\Theta$ and its inverse  $\Theta\in H\tens H$. We then eliminate the latter to give the conditions for a Hadamard form on $H$. We then write our conditions in terms of $\ch$ by inserting the red Frobenius form as needed and then recognise the dualisation of the red product and coproduct from Figure~\ref{figFalg} to obtain (2). \endproof
 
 \begin{example}\rm\label{exhZn} $H=\C\Z_n$ as an example of Example~\ref{exCX} has a quasi-Hadamard form as follows.  We write the group $X=\Z_n$ additively and to avoid confusion we let $|i\>$ explicitly denote the basis vector of $\C X$ spanned by $i\in \Z_n$ (a standard notation in the Computer Science literature). We let $q=e^{2\pi\imath\over n}$ be a primitive $n$-th root of unity and define
\[ \Theta(i,j)=q^{ij}\]
extended as a bilinear form on $H$. This is invertible and provides a Hadamard form. For example, on the part (a) side, $\Theta(i+j,k)=\Theta(i,k)\Theta(j,k)$ as required since $\Delta_{\color{green}\bullet}(|k\>)=|k\>\tens|k\>$. On the part (b) side, we have using $\Delta_{\color{red}\bullet}$,
\[\sum_{r+s=i}\Theta(s,j)\Theta(r,k)= \sum_rq^{(i-r)j}q^{rk}=\sum_rq^{r(k-j)}q^{ij}=n\delta_{j,k}q^{ij}=n\Theta(i,j{\color{green}\bullet}k)\]
where see the factor $n$ which makes this only true up to a constant factor. This is related to the fact that the red Frobenius algebra is quasi-special, not special. The associated Hadamard gate is 
\[ \ch(|i\>)=g^1_{\color{red}\bullet}\Theta(g^2_{\color{red}\bullet},|i\>)=\sum_j|j\> q^{-ji}\]
which is the usual $\Z_n$ (inverse) Fourier transform. Note for general $n$ that $\ch^2(|i\>)=n|-i\>$, so  $\ch$ can be normalised to be order $4$. For $n=2$ only, one can normalise by $1/\sqrt{2}$ to obtain $\ch=\ch^{-1}$, with both (a) and (b) conditions sharing the scale factor discrepancy.  \end{example}
 
Although this works better from the point of view of mapping red spiders to green spiders, it is not very natural from a Hopf algebra point to map both the algebra and coalgebra in the same direction. This suggests two natural variants.
 
 \begin{definition}\label{defself} Let $H$ be an $F$-Hopf algebra with antipode of the form in Proposition~\ref{propFhopf} and $\Theta$ a nondegenerate bilinear form on $H$. We call it a Type 2 (resp. Type 3) Hadamard form if it obeys (a) in Definition~\ref{defhad} and (b') (resp. (b'')), where
 \[ \includegraphics{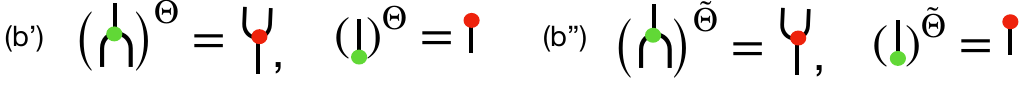}\]
and $\tilde\Theta=\Theta\circ(S\tens\id)$. 
 \end{definition}
 
Recall that when making a Hopf algebra $H$ into an $F$-Hopf algebra, we take it as the red product/green coproduct one and refer the green product/red coproduct one as the associated Hopf algebra $\tilde H$.

\begin{figure}
 \[ \includegraphics[scale=1]{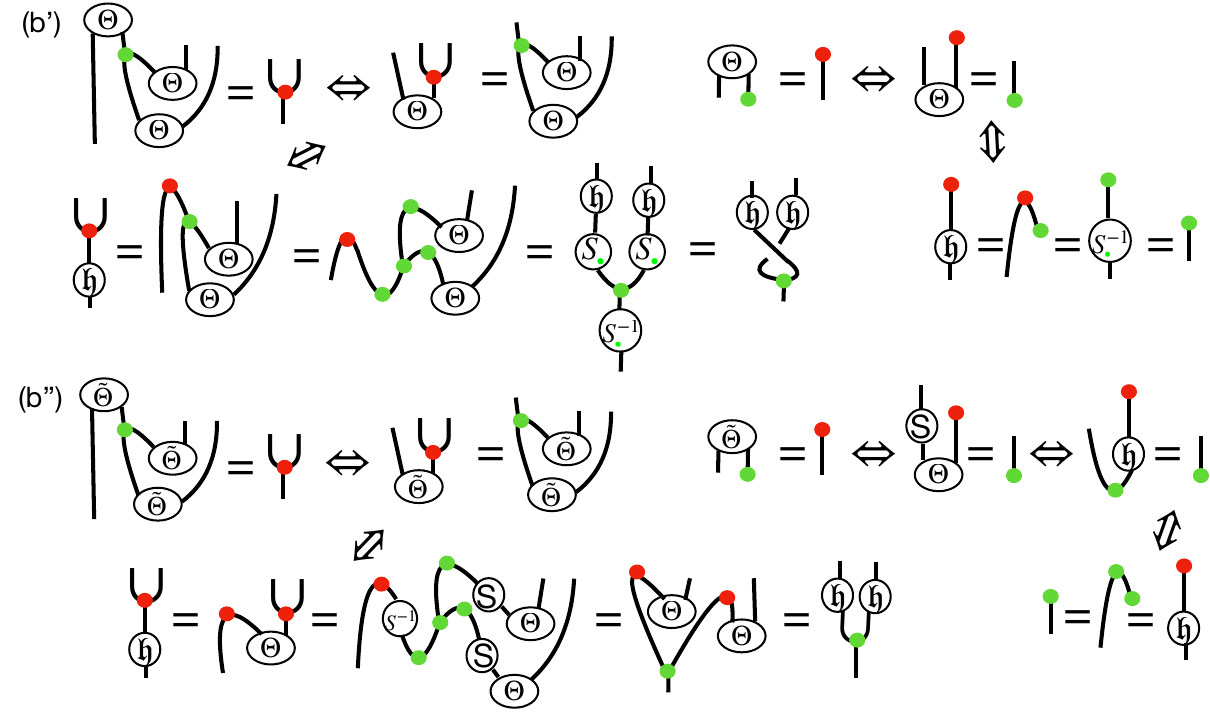}\]
 \caption{\label{fighadvarpf} Variants (b') and (b'') of a Hadamard form as Hopf algebra self-duality and anti-self-duality.}
 \end{figure}

\begin{lemma}\label{lemsd} Let $H$ be a Hopf algebra extended to an F-Hopf algebra by Proposition~\ref{propFhopf} and $\theta:H\tens H\to k$. The following are equivalent 

(1) $\Theta$ is a Type 2 (resp Type 3) Hadamard form.

(2) $h\mapsto \Theta(h, )$ is an isomorphism $H\isom H^\medstar$ (resp. $H\isom H^{\medstar op}$) of Hopf algebras.

(3) $\ch$ is and isomorphism  $H\cong \tilde H^{op}$ (resp. $H\isom \tilde H$) of Hopf algebras, i.e.
\[{\rm Type\ 2}:\quad \mu^{op}_{\color{green}\bullet}=\ch\circ\mu_{\color{red}\bullet}\circ(\ch^{-1}\tens\ch^{-1});\quad {\rm Type\ 3}:\quad \mu_{\color{green}\bullet}=\ch\circ\mu_{\color{red}\bullet}\circ(\ch^{-1}\tens\ch^{-1})\]
for the two cases and
\[  \Delta_{\color{green}\bullet}=(\ch^{-1}\tens\ch^{-1})\circ\Delta_{\color{red}\bullet}\circ\ch,\quad  1_{\color{green}\bullet}=\ch(1_{\color{red}\bullet}),\quad \eps_{\color{green}\bullet}=\eps_{\color{red}\bullet}\circ\ch.\]
\end{lemma}
\proof For Type 2, part (b') in Figure~\ref{fighadvarpf} shows the modified adjunction condition in (1) and converts it to a condition on $\Theta$. This version of conditions (a) and (b') explicitly is
\[ \Theta(hh',h'')=\Theta(h,h''_2 )\Theta(h',h''_1),\quad \Theta(h,h' h'')=\Theta(h_2,h')\Theta(h_1,h''),\quad \]
as well as $\Theta(1,h)=\eps(h)=\Theta(h,1)$, for all $h,h',h''\in H$. This is (2). In the figure, we then replace $\Theta$ by $\ch$ and recognise the antipode $\tilde S$ of the associated Hopf algebra $\tilde H$. We then use the (braided) anti-algebra homomorphism property of antipodes to obtain the variant in (3).  Similarly, for Type 3, part (b'') of Figure~\ref{fighadvarpf} shows the modified adjunction condition in (1) and converts it to a condition on $\tilde\Theta$. Then the (a) and (b'') conditions on the bilinear form say respectively that 
\[ \Theta(hh',h'')=\Theta(h,h''_2 )\Theta(h',h''_1),\quad \Theta(h,h' h'')=\Theta(h_1,h')\Theta(h_2,h''),\]
as well as $\Theta(1,h)=\eps(h)=\Theta(h,1)$, for all $h,h',h''\in H$. This is (2). In the figure, we then insert the assumed specific form of $S$ and recognise the green-dual of the green product of $H$ as needed for (3).  \endproof

Note that $H^\medstar\isom H^*$ as Hopf algebras (via the antipode), so Type 2 is equivalent to saying that $H$ is self-dual in the usual sense, while Type 3 is equivalent to saying that $H$ is anti-self-dual. 

  \begin{example}\rm $H=\C Z_n$ is a self-dual (or anti-self dual) Hopf algebra as the underlying finite Abelian group is self-dual. Here the condition (a) for  $\Theta(i,j)=q^{ij}$ is as in Example~\ref{exhZn} but instead of condition (b), we now have (b') which holds for the same reason as (a), given that the form is symmetric. Thus $\Theta$ is an exact (not quasi) Type 2  Hadamard form  leading more conventionally to $\ch$ as Fourier transform.  The actual $n=2$ Hadamard gate still changes the normalisation to have $\ch=\ch^{-1}$. \end{example}

\begin{example}\rm Let $H=u_q(b_+)$ be the Taft algebra with $q$ a primitive $n$-th root of unity. This was presented as an F-Hopf algebra in \cite{CoD}, and $u_q(b_+)$ is, moreover, known to be self-dual, see \cite[Ex~7.2.9]{Ma} at least for generic $q$. In our conventions, it is the sub-Hopf algebra of $u_q(sl_2)$ in Example~\ref{exsl2} generated by $K,F$ (say) and has
\[ \int K^iF^j=\delta_{i,0}\delta_{n-1,j},\quad \Lambda=\Lambda_KF^{n-1}.\]
We compute this for $n=3$. Then $\Lambda=(1+K+K^2)F^2$ and
\begin{align*} g_{\color{red}\bullet}&=F^2\tens 1+q^2K^2F\tens KF+q^2K\tens K^2F^2\\
&\quad +KF^2\tens K^2+F\tens F+qK^2\tens KF^2\\
&\quad  +K^2F^2\tens K+qKF\tens K^2F+1\tens F^2\end{align*}
\[  (K^iF^m){\color{green}\bullet}(K^jF^k)=\delta_{[i-j],2-k}\left[{m\atop 2-k}\right]_q(-1)^kq^{\delta_{k,1}}K^jF^{m+k-2}\]
treated as 0 unless $m+k\ge 2$. Here $i,j$ are mod 3 and $[i-j]$ denotes the value in the range 0,1,2 and the $q$-binomials are defined using $q$-integers. We used the $q$-binomial formula to compute $\Delta(K^iF^m)$ and pick off the term with the right power of $F$ for the integral to not vanish. We let 
\[ t=(1+qK +q^2K^2)F^2,\quad x=(1+K+K^2)F\]
then we find
\[ t{\color{green}\bullet}t{\color{green}\bullet}t=\Lambda,\quad x{\color{green}\bullet}t=(1+qK+q^2K^2)F,\quad t{\color{green}\bullet}x=qx{\color{green}\bullet}t,\quad x{\color{green}\bullet}x{\color{green}\bullet}x=0,\]
where $\Lambda=1_{\color{green}\bullet}$. The red coproduct is given by multiplication by $g_{\color{red}\bullet}$ but this is central, so we can conveniently compute it as
\[\Delta_{\color{red}\bullet}(K^iF^j)=K^ig_{\color{red}\bullet}F^j;\quad \Delta_{\color{red}\bullet} t=t\tens t,\quad \Delta_{\color{red}\bullet}x=x\tens t^{-1}+ \Lambda\tens x.\]
 Thus, $H$ is isomorphic to the opposite of the associated Hopf algebra by  $K\mapsto t^{-1}$ and $F\mapsto t{\color{green}\bullet} x$, so we have a Type 2 Hadamard form. We can now read off some of the values of $\ch$ and deduce the rest from the algebra homomorphism property. Thus, on the vector space of $H$, 
 \begin{align*} \ch(K^i)&=\delta_{-i}F^2,\quad
 \ch(K^iF)=q\delta_{1-i} F,\quad   \ch(K^iF^2)=\delta_{2-i},\end{align*}
where we define
\[ \delta_i=\sum_j q^{ij} K^j.\]
We see some similarities to the $\Z_n$ Fourier transform.  \end{example}

\begin{example}\rm Let $H=u_{-1}(sl_2)$ as studied in Example~\ref{exsl2}. This is both self dual and anti-self dual as a Hopf algebra. Here, we already gave an isomorphism from $H$ to the associated Hopf algebra by specifying it on the generators, and deduce the rest now as
\[ \ch(1)=(1+K)EF,\quad \ch(K)=(K-1)EF,\quad\ch(E)=(1+K)E,\quad\ch(F)=(1-K)F\]
\[ \ch(KE)=(K-1)E,\quad \ch(KF)=(1+K)F,\quad \ch(FE)=-K,\quad \ch(KFE)=1+K\]
as the corresponding Type 3 Hadamard gate. 
\end{example}

The last example here as well as $\C Z_n$ are flip-Hopf *-algebras, so it is reasonable to ask that $\Theta$ and  $\ch$ be compatible with the $*$-algebra structure when they exist. We do not have a general analysis for this, but for $\C \Z_n$ we have $\Theta(i,j)=q^{ij}$ as a matrix is unitary up to a constant factor $n$ in our orthonormal basis defined from the $*$-structure. It is not, however, the case in this example that $\ch$ as an operator is a $*$-algebra map, but that is the case in the $u_{-1}(sl_2)$ example. 

To close this section, we note that for a more categorical view on Frobenius $*$-algebras and F-Hopf $*$-algebras, one cannot work in the category of vector spaces alone since $*$ is not a morphism of vector spaces. An approach that is useful in noncommutative geometry is the notion of a bar-category\cite{BegMa:bar}\cite[Chap.~2.8]{BegMa}, in which a monoidal category is equipped with a functor `bar' that conjugates objects. Here every complex vector space $V$ has a conjugate $\bar V$ which is the same abelian group but $z\in \C$ acts by its conjugate. Then for a $*$-algebra $A$,  $*:A\to \bar A$ is a linear map. There are related notions in the Computer Science literature. 

\section{Braided interacting Hopf algebras}\label{secbraZX}

Braided ZX calculus was considered in \cite{BPW}, where it was shown that the bialgebra rewrite rules as stated there force the category to be symmetric. Here we revisit this from the point of view of building up the required algebraic axioms or `rewrite rules' by understanding the correct construction of a braided interacting pair of Hopf algebras. 

We recall that a braided category is first of all a monoidal category, so there is a functor for $\tens$ of objects and morphisms, with unit object $\underline 1$ and associated maps including an natural equivalence $(\ \tens\ )\tens \to \tens(\ \tens\ )$ (the associator). The Mac Lane coherence theorem for monoidal categories says that we can ignore brackets and insert the associator as needed.  In addition, there is a natural equivalence $\Psi: \tens\to \tens^{op}$ (the braiding). The coherence theorem for braided categories was in \cite{JoyStr}, replacing symmetric group by the braid group in Mac Lane's coherence theorem for symmetric monoidal categories. This underlies the diagrammatatic notation for working in braided categories, with $\Psi=\includegraphics{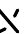}$ a braid crossing, which in turn underlies the theory of braided-Hopf algebras\cite{Ma:bra,Ma:tra,Ma:cro,Ma:alg}. The axioms for these are as in Section~\ref{secFHopf}, but now caring about over- and under- crossings. The most important thing we need from the theory of braided-Hopf algebras is a result\cite{Ma:pri,Ma:alg} that 
\begin{equation}\label{antialg} \includegraphics{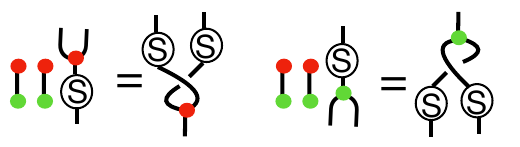}\end{equation}
where we have inserted the factor for an unnormalised Hopf algebra. Here $S$ denotes the antipode of a braided Hopf algebra $H$ with red product and green coproduct, not necessarily of any special form. We proceed in the normalised case of a usual braided Hopf algebra. 
We already used this {\em braided anti-homomorphism} property in the proofs in the preceding section, with the braided theory of the present section in mind. 

The definition of an F-Hopf algebra still makes sense declaring $H$ to be braided, but we saw in the proof of Proposition~\ref{propFhopf} {\em one spot} where we needed a wrong braid crossing.  Thus, what we actually proved in parts (a)-(b) of Figure~\ref{figFhopf} was the following.  

\begin{corollary} Let $H$ be a Hopf algebra in a braided category $\CC$ with red algebra and green coalgebra Frobenius and with antipode of the form in Proposition~\ref{propFhopf}. Then the associated green algebra and red coalgebra form a Hopf algebra $\tilde H$ in the braided category with reversed braid crossings (or a reverse-braided Hopf algebra in $\CC$). 
\end{corollary}

In fact, the opposite algebra to a Hopf algebra with reversed braiding is a braided Hopf algebra in the original category\cite{Ma:pri}. Equivalently,  the categorical dual $H^\medstar$ of $H$ remains in $\CC$\cite{Ma:alg,Ma:pri} and part (c) of Figure~\ref{figFhopf} says that if we compute this using the red duality then we indeed land on the opposite product. Hence, we also proved the following.

\begin{corollary}\label{propbraFhopf} If a braided Hopf algebra $H$ is Frobenius as an algebra and coalgebra and its antipode $S$ has the special form in Proposition~\ref{propFhopf} then 
\[ \includegraphics[scale=0.85]{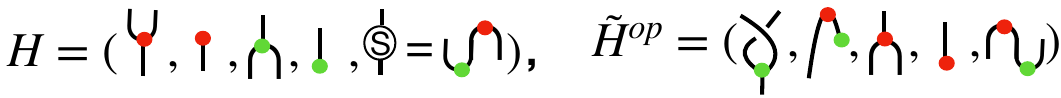}\]
are a pair of braided Hopf algebras on the same object $H$ while the red and green maps are a pair of F-algebras. The second braided Hopf algebra is isomorphic to $H^\medstar$. \end{corollary}

Either point of view gives a canonical braided F-Hopf algebra amplified under our assumptions from a single braided Hopf algebra. Note that, as with ordinary F-Hopf algebras, it is not clear that there are many examples beyond  these,  so this is more of a construction rather than a general concept of interacting Hopf algebras.  We let $\underline 1$ be the unit object of the braided category. 

\begin{figure}
 \[ \includegraphics[scale=1]{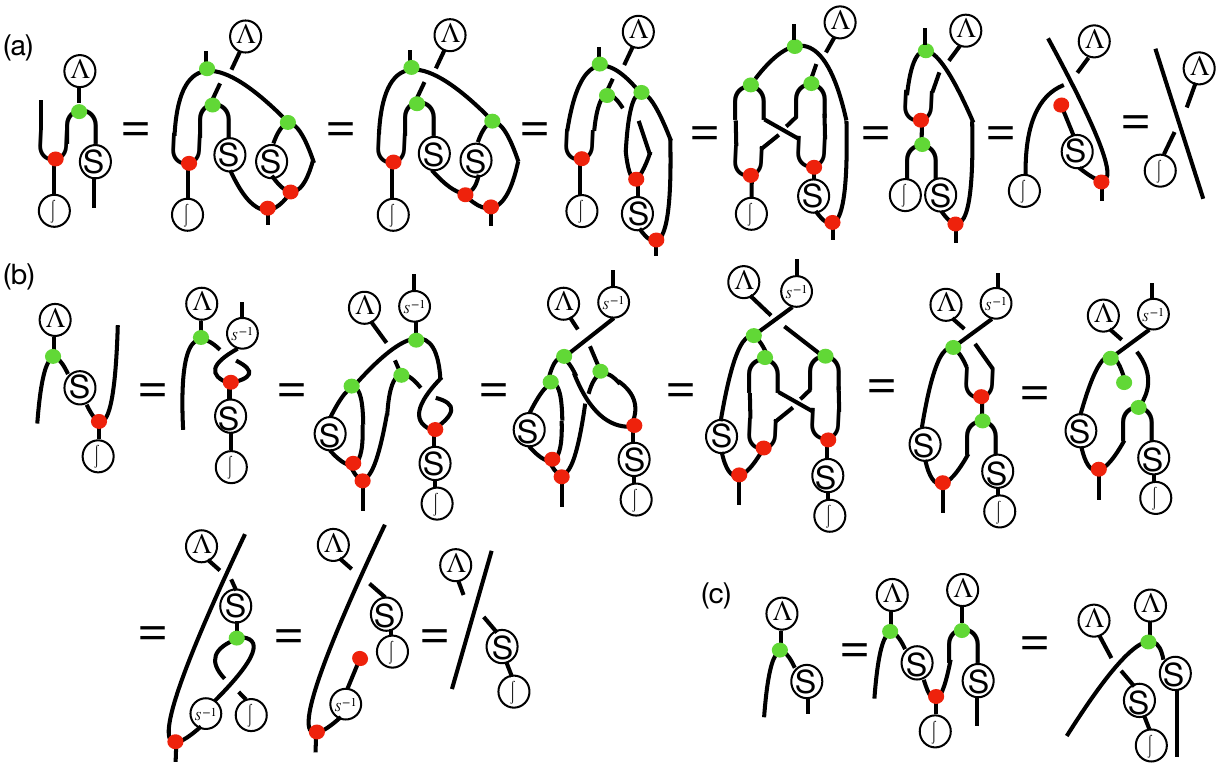}\]
 \caption{\label{figbraint} Construction of Frobenius form from  integrals $\Lambda,\int$ on a braided Hopf algebra.}
 \end{figure}

\begin{corollary}\label{corbraHH*} Let $H$ be a braided-Hopf algebra with invertible antipode equipped with a morphism $\Lambda:\underline 1\to H$ which is a left integral  and a morphism $\int:H\to \underline 1$ which is a right integral, such that $\int\Lambda$ is the identity on $\underline 1$. Then  Proposition~\ref{propbraFhopf} applies in the braided category and we have a braided F-Hopf algebra.
\end{corollary}
\proof We follow the definitions and steps in the proof of Corollary~\ref{corHH*} but now as a morphism in a braided category. Thus, we set
\[ (\ ,\ )_{\color{red}\bullet}=\int\circ\mu_{\color{red}\bullet},\quad  g_{\color{red}\bullet}=(\id\tens S)\Delta\Lambda;\quad  (\ ,\ )_{\color{green}\bullet}=(\ ,\ )_{\color{red}\bullet}(S\tens\id),\quad g_{\color{green}\bullet}=(\id\tens S^{-1})g_{\color{red}\bullet}.\]
The proof for the red Frobenius form is then shown in Figure~\ref{figbraint}. Part (a) uses the braided-antimultiplicativity property for the 3rd equality and the integral property of $\int$ for the 6th. At the end, we need that $\Lambda$ or $\int$ are morphisms to obtain the identity map. Part (b) uses braided-antimultiplicativity for the 1st equality, the morphism and integral property of $\Lambda$ for the 6th and the morphism and integral property of $\int$ for the 8th. These results together imply in (c) that $\int S\Lambda$ is the identity on $\underline 1$ under either morphism assumption. It is immediate that we have the duality properties for the green form and that this makes the coalgebra Frobenius, and that $S$ then has the form required in Proposition~\ref{propbraFhopf}.\endproof

In fact, integrals on braided-Hopf algebras are not typically bosonic (i.e. of trivial braiding), so this involves a certain complication already familiar in the literature on braided Fourier transform\cite{BesLyu}. This remark already applies to the {\em braided line} $\C[x]/(x^n)$ with degree $|x|=-1$ in the category of $\Z_n$-graded spaces with primitive coalgebra and antipode
\[ \Delta x^m=\sum_{r=0}^m\left[{m\atop r}\right]_q x^r\tens x^{m-r},\quad S x^m= q^{m(m-1)\over 2}(-x)^m,\quad \eps x^m=\delta_{m,0}\]
extended with braiding $\Psi(x^i\tens x^j)=q^{ij}x^j\tens x^i$. This has an integral $\int: H\to K$ where $K=k$ regarded as a 1-dimensional object of grade -1 and given by $\int x^m=\delta_{n-1,m}$ and an integral element $\Lambda=x^{n-1}$ also of grade -1. As a result, the corollary does not apply on this simplest example, i.e. 
\[ (x^a,x^b)_{\color{red}\bullet}=\delta_{n-1,a+b},\quad g_{\color{red}\bullet}=\sum_a x^a\tens x^{n-1-a},\quad \mu_{\color{red}\bullet}(g_{\color{red}\bullet})=n x^{n-1}\]
where $g_{\color{red}\bullet}$ provides the inverse to $(\ ,\ )_{\color{red}\bullet}$ but is not of the expected form $(\id\tens S)\Delta x^{n-1}$. Likewise,
\[ (x^a,x^b)_{\color{green}\bullet}={\delta_{n-1,a+b}\over \left[{n-1\atop a}\right]_q},\quad g_{\color{green}\bullet}=\Delta x^{n-1}=\sum_{a=0}^{n-1}\left[{n-1\atop a}\right]_q x^a\tens x^{n-1-a}\]
where $(\ ,\ )_{\color{green}\bullet}$ provides the correct inverse to $g_{\color{green}\bullet}$, but is not of the expected form $(S(\ ),\ )_{\color{red}\bullet}$. As a result, $S$ is not of the required form for the braided reading of Proposition~\ref{propbraFhopf}.

\subsection{Self-dual braided interacting pairs by transmutation} \label{sectrans}

On the other hand, all quasitriangular Hopf algebras, such as $u_q(sl_2)$, have a self-dual braided Hopf algebra version via a process of {\em transmutation}\cite{Ma:bra,Ma:tra,Ma,Ma:pri}. We show that these provide examples for the braided theory. They also underly TQFT invariants and hence should  be relevant to fault tolerant TQFT based quantum computing. 

Let $(H,\CR)$ be an ordinary  quasitriangular Hopf algebra\cite{Dri,Ma}. Its transmutation $\underline H$ has the same algebra and counit but a different coproduct and antipode\cite{Ma:bra,Ma:tra}
\begin{equation}\label{trans} \underline{\Delta}h=h_1S\CR^2\tens \CR^1\la h_2, 
\quad \underline S h=\CR^2 S(\CR^1\la h)\end{equation}
forming a braided-Hopf algebra in the category of $H$-modules as an object by the adjoint action
\begin{equation}\label{adj} h\la h'= h_1 h' S h_2.\end{equation}
Throughout this section, we will use underlines to separate the braided Hopf algebra structures from  the original ordinary ones. We will also be interested in the {\em factorisable} case which, by definition, means that the quantum Killing form $\CQ:=\CR_{21}\CR$ is invertible when viewed as a linear map $H^*\to H$ by evaluation against its first factor. Here, if $S\Lambda=\Lambda$ then $h\la \Lambda=\eps(h)\Lambda$, so we can regard it as a morphism from the trivial object to $\underline{H}$. The following is a reworking of arguments in \cite{LyuMa}. 

\begin{proposition}\label{proptransH} Let $H$ be a factorisable quasitriangular Hopf algebra with $S\Lambda=\Lambda$. Then Corollary~\ref{corbraHH*} applies to $\underline H$ and provides a braided F-Hopf algebra in the braided category of $H$-modules. 
\end{proposition}
\proof  If $S\Lambda=\Lambda$ then $h\la \Lambda=\eps(h)\Lambda$, so we can regard it as a morphism from the trivial object to $\underline{H}$. Next, $H^*$ is coquasitriangular and a dual transmutation construction keeps the same coalgebra but modifies the algebra, making it into a braided-Hopf algebra, which can be viewed in the same braided category by the quantum coadjoint action. In this case, $\CQ$ as a map becomes an isomorphism of braided-Hopf algebras, see \cite{Ma,Ma:pri} for more details. Finally, the right integral $\int$ on $H^*$ is also a right integral on the braided version and therefore transfers under the braided Hopf algebra self-duality isomorphism to a  right  integral $\underline\int$ on $\underline H$. By arguments \cite[Lem.~3.5]{LyuMa}, it is ad-invariant and hence a morphism. Therefore, we meet the conditions of Corollary~\ref{corbraHH*}. \endproof

Rather than more formal details, will now show how this all works in the simplest quantum group example. 

\begin{example}\rm Let $q$ be a primitive $n$-th root of unity {\em with $n$ odd} and $H=u_q(sl_2)$ as in Example~\ref{exsl2}. In this case $q$ has a square root $q^{n+1\over 2}$ and $u_q(sl_2)$ in our conventions (denoted by gothic in \cite{AziMa}) is isomorphic to $u_{q^{1\over 2}}(sl_2)$ in usual conventions as shown there. Likewise for our conventions for the reduced $c_q[SL_2]$ as its Hopf algebra dual. The braided Hopf algebra $\underline{H}=b_q(sl_2)$ has the same algebra as $u_q(sl_2)$ and hence $\Lambda=\delta_0(K)F^nE^n$ there is also a left integral in $H$, where we use the notation $\delta_i:=\sum_i q^iK^i$. We have $S\Lambda=\Lambda$, so this is invariant under the adjoint action, which is equivalent to $\Lambda\in Z(H)$. The latter can also be checked from the commutation relations. The coalgebra of $b_q(sl_2)$ is obtained from that of $u_q(sl_2)$ via (\ref{trans}) using the adjoint action and the quasitriangular structure. 

Also when $n$ is odd,  $u_q(sl_2)$ is factorisable and we have a self-duality braided-Hopf algebra isomorphism $b_q(sl_2)\isom b_q[SL_2]$, where $b_q[SL_2]$ is obtained by a dual transmutation construction from. It is built on the same coalgebra and viewing the matrix of generators of $c_q[SL_2]$ in $b_q[SL_2]$, where we denote them $\alpha,\beta,\gamma,\delta$, these have the relations, coproduct and antipode
\[ \alpha^n=1,\quad\beta^n=\gamma^n=0;\quad\quad\beta\alpha=q\alpha\beta,\quad \gamma\alpha=q^{-1}\alpha\gamma,\quad\delta\alpha=\alpha\delta, \quad \alpha\delta-q\gamma\beta=1\]
\[  [\beta,\gamma]=(1-q^{-1})\alpha(\delta-\alpha),\quad [\gamma,\delta]=(1-q^{-1})\gamma\alpha,\quad [\delta,\beta]=(1-q^{-1})\alpha\beta, \]
\[\underline{\Delta}\begin{pmatrix}\alpha&\beta\\ \gamma & \delta\end{pmatrix}=\begin{pmatrix}\alpha&\beta\\ \gamma & \delta\end{pmatrix}\tens \begin{pmatrix}\alpha&\beta\\ \gamma & \delta\end{pmatrix},\quad \underline{S}\begin{pmatrix}\alpha&\beta\\ \gamma & \delta\end{pmatrix}=\begin{pmatrix}q\delta+(1-q)\alpha&-q\beta\\ -q\gamma & \alpha\end{pmatrix}\]
and $\underline\eps(\alpha)=\underline\eps(\delta)=1, \underline\eps(\beta)=\underline\eps(\gamma)=0$. We replaced $q^2$ in \cite{Ma} by $q$ to match our conventions. The braiding $\Psi$ on the matrix of generators is given explicitly in \cite[Ex~10.3.3]{Ma}, upon making the same change from $q^2$ to $q$. Here
\[ \Psi(\alpha\tens\alpha)=\alpha\tens\alpha+(1-q)\beta\tens\gamma,\quad \Psi(\alpha\tens\beta)=\beta\tens\alpha,\quad\Psi(\alpha\tens\gamma)=\gamma\tens\alpha+(1-q)(\delta-\alpha)\tens\gamma\]
\[ \Psi(\beta\tens\alpha)=\alpha\tens\beta+(1-q)\beta\tens(\delta-\alpha),\quad\Psi(\beta\tens\beta)=q\beta\tens\beta,\]\[\quad \Psi(\beta\tens\gamma)=q^{-1}\gamma\tens\beta+(1+q)(1-q^{-1})^2\beta\tens\gamma-(1-q^{-1})(\delta-\alpha)\tens(\delta-\alpha)\]
\[ \Psi(\gamma\tens\alpha)=\alpha\tens\gamma,\quad\Psi(\gamma\tens\beta)=q^{-1}\beta\tens\gamma,\quad\Psi(\gamma\tens\gamma)=q\gamma\tens\gamma\]
with the braiding on products determined by iterated $\Psi$ as the algebra product is a morphism and hence passes through braid crossings. For example, $\Psi(\gamma\tens\delta)=\delta\tens\gamma$ and $\Psi(\delta\tens\beta)=\beta\tens\delta$ can be deduced.

The self-duality isomorphism to $b_q(sl_2)$ is likewise given in our conventions by
\[ \begin{pmatrix}\alpha&\beta\\ \gamma & \delta\end{pmatrix}\mapsto \begin{pmatrix}K& (1-q^{-1})E\\
(q^{-2}-1)KF & K^{-1}+(1-q)(1-q^{-2})FE\end{pmatrix}.\]
Moreover, as $b_q[SL_2]$ has the same coalgebra, it gets the right-invariant integral of $c_q[SL_2]$. Hence under the isomorphism, $b_q(sl_2)$ gets a right-invariant integral $\underline\int$ at least for odd $n$. Also note that the element $c$ defined as the image under the isomorphism of  $\alpha+q\delta$ (a multiple of the $q$-trace), namely
\[ c=K+qK^{-1}+(1-q)(q-q^{-1})FE\]
is central and hence invariant. Thus, $1,c,c^2,\Lambda$ are bosonic (they have trivial braiding with everything) and we also necessarily have $\underline Sc=Sc=c$ and $\underline S\Lambda= S\Lambda=\Lambda$ using the general formula for $\underline S$. 

Finally, we start with $\int$ on $c_q[SL_2]$ as in \cite{AziMa}, which when converted to the matrix generators $a,b,c$ (the fourth matrix generator $d$ is determined in the root of unity case) amounts to
\[ \int a^ib^jc^k=\delta_{i,0}\delta_{n-1,j}\delta_{n-1,k}.\]
We take take this  $b_q[SL_2]$ also as it has the same vector space and map it over to $b_q(sl_2)$ under the isomorphism of braided Hopf algebras. Corollary~\ref{corbraHH*} then applies and we have a Frobenius structures for the algebra and coalgebra to give a braided F-Hopf algebra. 

In what follows, we focus on  $n=3$, so $q$ is a primitive cube root of unity and $1+q+q^2=0$. Using the transmutation formulae \cite[Eqn (7.37)]{Ma} extended to degree 4 gives $c^2b^2=q c\bullet c\bullet b\bullet b=q\gamma^2\beta^2$ where the modified product is denoted $\bullet$ and we then view the elements in $b_q[SL_2]$. Thus, the braided integral is zero except on $\gamma^2\beta^2$ in the monomial basis for $\gamma,\beta,\alpha$. Normalised so that $\underline\int \Lambda=1$, we have after the isomorphism with $b_q(sl_2)$,
\[ \underline\int K^iF^jE^j=\delta_{2,i}\delta_{2,j}\delta_{2,k}\]
as the braided integral. It is easy to see that this is invariant in a Hopf algebra sense under the quantum adjoint action and hence a morphism. Also note that this is not the same as the integral on $u_q(sl_2)$ as the coproduct is different, but similar. 

We will now do some checks of the construction. From the matrix form of braided coproduct of $b_q[SL_2]$ and the stated isomorphism, we can read off some coproducts and antipodes
\begin{align*}
\underline\Delta K&=K\tens K+q^{-1}(q-1)^2E\tens KF,\quad \underline S K=c-qK,\\
 \underline\Delta E&=E\tens q^{-1}(c-K)+K\tens E,\qquad \underline S E=-qE,\\
\underline\Delta(KF)&=KF\tens K+q^{-1}(c-K)\tens KF,\quad\underline S(KF)=-qKF\end{align*}
where $\delta$ maps under the isomorphism to $q^{-1}(c-K)$. 
The second of these is relatively easy to compute directly using the first four of 
\begin{gather*}K\la K=K,\quad K\la(KF)=q KF,\quad \\ F\la K=(q^{-1}-1)KF,\quad F\la(KF)=0,\quad F\la K^2=(q-1)K^2F\end{gather*}
and $\CR$ from Example~\ref{exsl2}, as higher order terms of $\CR$ do not contribute. One can also check that the algebras are isomorphic. Some examples of the braiding can likewise be read off as
\begin{align*} \Psi(K\tens E)&=E\tens K,\quad \Psi(E\tens E)=qE\tens E,\quad \Psi(KF\tens K)=K\tens KF,\\
\Psi(KF\tens E)&=q^{-1}E\tens KF,\quad \Psi(KF\tens KF)=qKF\tens KF,\\
 \Psi(E\tens K)&=K\tens E+ (1-q)E\tens (q^{-1}c+qK),\\
 \Psi(K\tens KF)&=KF\tens K+(1-q)(q^{-1}c+qK)\tens KF\\
  \Psi(K\tens K)&=K\tens K-q^{-1}(q-1)^3E\tens KF,\\
  \Psi(E\tens KF)&=q^{-1}KF\tens E-(q-1)^2E\tens KF-{1\over q-1}(q^{-1}c+qK)\tens (q^{-1}c+qK)
 \end{align*}
 which allows us to compute more coproducts  and antipodes using the braided-homomorphism and anti-homomorphism properties respectively. For example,
\begin{align*}\underline S(E^2)&=\cdot\Psi((-qE)\tens(-qE))=E^2\\
\underline S (K^2F^2)&=q\cdot\Psi( (-q KF)\tens(-qKF))=qKFKF=K^2F^2\\
\underline S(KE)&=\cdot\Psi((c-qK)\tens(-qE))=qE(qK-c)=(K-qc)E\\
 \underline S (K^2)&=\cdot\Psi((c-qK)\tens(c-qK))= c^2+q^2(K^2+(1-q)(1-q^{-1})(q^{-2}-1)EKF)- 2 q c K\\
&=c^2+ c K + K^2 + q(q-1)=:s(c,K)\\
\underline S(K^2F)&=q\underline S((KF)K)=q\cdot\Psi((-qKF)\tens(c-qK))=-q^2(c-qK)KF=(K-q^2c)KF\\
\underline S(KE^2)&=\cdot\Psi((c-qK)\tens E^2)=E^2(c-qK)=(c-K)E^2\\
\underline S(F^2)&=\underline S((KFKF)K)=\cdot\Psi(KFKF\tens(c-qK))=(c-qK)KFKF=(q^{-1}cK^2-1)F^2\\
\underline S F&=q^{-1}\underline S(KFK^2)=q^{-1}\cdot\Psi((-qKF)\tens s(c,K))=- s(c,K)KF\\
\underline S(KFE)&=\underline S((KF)E)=\cdot \Psi((-qKF)\tens(-qE))=q EKF=q^2 KEF\\
\underline S(F^2E^2)&=q^2\underline S(K^2F^2KE^2)=\cdot(KFKF\tens (c-K)E^2)=q^{-1}(c-K)E^2KFKF=(cK^2-1)E^2F^2\\
\underline S(K^2E)&=\cdot\Psi(s(c,K)\tens (-qE))=-q Es(c,K)\\
\underline S(K^2E^2)&=\cdot\Psi(s(c,K)\tens E^2)=E^2s(c,K)\\
\underline S(FE^2)&=q^{-1}\underline S(KF(K^2E^2))=-\cdot\Psi(KF\tens E^2 s(c,K))=-q^{-2}E^2  s(c,K) KF\\
\underline S(KF^2E)&=q^{-1}\underline S(KFKF(K^2 E))=-\cdot\Psi(KFKF\tens Es(c,K))=-q^{-2}E s(c,K)KFKF\\
&=-E s(c,K)K^2F^2.
\end{align*}

We now let $(\ ,\ )_{\color{red}\bullet}=\underline\int\mu_{\color{red}\bullet}$ and $g_{\color{red}\bullet}=(\id\tens\underline S)\underline\Delta \Lambda$. By construction, these will be inverse to each other and
as a spot check (to make sure there was no confusion in our conventions) we demonstrate this on a sample element, $F$, say, i.e. we check that
$(F,g^1_{\color{red}\bullet})_{\color{red}\bullet} g^2_{\color{red}\bullet}=F$. The same method will apply on other elements similarly. We start with  the unbraided coproduct
\begin{align*} \Delta\Lambda&=(1\tens 1+K\tens K+K^2\tens K^2)(F^2\tens 1+K\tens F^2-q^2K^2F\tens F)\\
&\qquad\qquad\qquad (1\tens E^2+E^2\tens K^2-q^2E\tens KE)\\
&=(1\tens 1+K\tens K+K^2\tens K^2)\big(F^2\tens E^2+F^2E^2\tens K^2-q^2F^2E\tens KE\\
&\qquad\qquad\qquad+K\tens F^2E^2+q^2KE^2\tens K^2F^2-KE\tens KF^2E\\ &\qquad\qquad\qquad -q^2K^2F\tens FE^2-K^2FE^2\tens K^2F+K^2FE\tens KFE\big)\end{align*}
and we now compute
 \begin{align*}\Big( \underline\int F&\Lambda_1 S \CR^2 \Big) \underline S(\CR^1\la \Lambda_2)=\left({1\over 3}\sum_{a,b,r}\int {(-1)^r(q-q^{-1})^r\over[r]_{q^{-1}}!} q^{-ab} F\Lambda_1 S(E^r K^b)\right) \underline S((F^r K^a)\la \Lambda_2)\\
 &=\left(\sum_{r}\int {(q-q^{-1})^r\over[r]_{q^{-1}}!} F\Lambda_1 K^2q^{-{r(r+1)\over 2}}E^r\right) \underline S(F^r\la \Lambda_2)\\
 &=\left(\int F(-FE^2)K^2\right)\underline S F+q^{-1}(q-q^{-1})\left(\int F(FE)K^2E\right)\underline S F\la(K^2FE)\\
 &\qquad+{(q-q^{-1})^2q^{-3}\over1+q^{-1}}\left(\int F( -q^2F)K^2E^2\right)\underline S F^2\la(KFE^2)\\
 & = \underline S\big( -F+q(1-q)F\la(K^2FE)+(1-q)^2F^2\la(KFE^2)\big)\\
 &=\underline S\big( -F+(1-q)^{-1}F\la(c K^2-1-qK)+F^2\la((q^2 cK-q^2 K^2-1)E)\big)\\
 &=\underline S\big( -F+(1-q)^{-1}F\la(c K^2-q K-q c^2 K+q^2K)\big)\\
 &=\underline S( -F+ (q-1-c^2-cK)KF)=F\end{align*}
 as required. Here, in order to have a nonzero integral, $\Lambda_1$ must involve $F E^{2-r}$, which means $\Lambda_2$ has $F E^r$ and hence $K^a\la \Lambda_2=q^{(1-r)a}\Lambda_2$. We then do the sum over $a$ which absorbs the 1/3 and sets $b=1-r$, which is our second expression on using $S E^r=(-1)^r q^{-{r(r+1)\over 2}}K^{-r}E^r$. For the example of $F$ that we are testing on, the relevant terms in the 2nd displayed factor of $\Delta\Lambda$ are one of $ -K^2FE^2, K^2FE,-q^2K^2F$ according to $r=0,1,2$. In each case the contributing from the prefactor in $\Delta\Lambda$ must be $K\tens K$ in order for the total power of $K$ in the integrand to be 2. This gives the 3rd expression. Evaluating the integrals gives the 4th expression and we then replace $FE$ in terms of $c$ for the 5th. Next,
\begin{align*}F\la &((cK-K^2-q)E)=F(cK-K^2-q)E+K^{-1}(cK-K^2-q)E(-KF)\\
&=(cq^{-1}-qK^2-q)FE-q(cK-K^2-q)EF={1\over q-1}((c^2-q)K-q c-1)\end{align*}
which leads to the 6th expression when the further  $F\la$  is combined with the existing one. We then compute $F\la$ noting that $c$ is invariant and $F\la K, F\la K^2$ were already given above. We then compute $\underline S$, noting that $c$ is unchanged and passes through it, and using $\underline S F, \underline S(KF)$ and $\underline S(K^2F)$ already found above. In this way, we obtain $F$ as expected. 

Replacing $F$ by a different monomial will pick of 3 different terms for $\Lambda_1$, etc, i.e. we can proceed similarly. As in the proof of Corollary~\ref{corbraHH*}, the green Frobenius structures  follow from the red ones via $\underline S$ and its inverse, to complete the construction. Moreover, our arguments apply for general odd $n$.   \end{example}

Proposition~\ref{proptransH} also applies in principle to other reduced `Lusztig kernel' quantum groups $u_q(\cg)$, although the precise version of these needed and at which roots of unity is not fully understood, see \cite{AziMa}. Returning to the general construction of $\underline H$ in the factorisable case, as these are self-dual, we also have a Type 2 Hadamard form in the sense of Lemma~\ref{lemsd}. 

\begin{proposition} Let $(H,\CR)$ be a factorisable quasitriangular Hopf algebra. Then there is a Type 2 Hadamard form $\underline H\tens \underline H\to \underline 1$ according to the braided version of Lemma~\ref{lemsd}, inverse to the  metric element
\[ \Theta=(S\tens\id)\CQ\in \underline H\tens\underline H.\]
\end{proposition}
\proof This is a version of the self-duality isomorphism between the the transmutation of $H$ and the transmutation of $H^*$, which was used in the construction above. It is known \cite{LyuMa,Ma} that $(S\tens\id)\CQ$ is  quantum ad-invariant and hence a morphism from $\underline 1$. That this is invertible when viewed as a map $\underline H^\medstar\to\underline H$ is essentially the assumption of being factorisable. The diagrams in Figure~\ref{fighadpf}(a) and Figure~\ref{fighadvarpf} part (b') when written in terms of the metric element $\Theta$ (rather than the bilinear form itself) amount to 
\[ (\underline\Delta\tens\id)\Theta=(\id\tens\mu)\Theta_{23}\Theta_{13},\quad (\id\tens\underline\Delta)\Theta=(\mu\tens\id)\Theta_{12}\Theta_{23},\quad(\id\tens\eps)\Theta=(\eps\tens\id)\Theta=1\]
The counit parts are immediate from $(\eps\tens\id)\CR=(\id\tens\id)\CR=1$ while 
\begin{align*} (\underline\Delta&\tens\id)\Theta=S(\CR^{iv 2}\CR^2{}_2\CR^{'1}{}_2)\tens \CR^{iv 1}{}_1S(\CR^{iv 1}{}_2 \CR^2{}_1\CR^{'1}{}_1)\tens\CR^1\CR^{'2}\\
&=S(\CR^{iv2}\CR^{v2}\CR^{''2}\CR^{'''1})\tens \CR^{iv 1} S(\CR^{v 1} \CR^2\CR^{'1})\tens\CR^{''1}\CR^1\CR^{'2}\CR^{'''2}\\
&=S(\CR^{iv2}\CR^{''2}\CR^{v2}\CR^{'''1})\tens \CR^{iv 1} S( \CR^2\CR^{v 1}\CR^{'1})\tens\CR^1\CR^{''1}\CR^{'2}\CR^{'''2}\\
&=S(\CR^{iv2}\CR^{v2}\CR^{''2}\CR^{'''1})\tens \CR^{iv 1} S( \CR^2\CR^{'1}\CR^{v 1})\tens\CR^1\CR^{'2}\CR^{''1}\CR^{'''2}\\
&=S(\CR^{''2}\CR^{'''1})\tens S( \CR^2\CR^{'1})\tens\CR^1\CR^{'2}\CR^{''1}\CR^{'''2}=\Theta_{23}\Theta_{13}\\
\end{align*}
where $\CR=\CR^1\tens\CR^2$ (summation understood) and the primed/roman superfixes denote independent copies. We used (\ref{trans}) for $\underline\Delta$ and properties of the antipode for the first equality, the quasitriangularity properties  $(\Delta\tens\id)\CR=\CR_{13}\CR_{23}$ and $(\id\tens\Delta)\CR=\CR_{13}\CR_{12}$ for the 2nd, then the braid relations for $\CR$ for the 3rd and 4th. For the 5th we use $(S\tens\id)\CR=\CR^{-1}$ applied to $\CR^{v}$ and cancelled with $\CR^{iv}$. Similarly on the other side,
\begin{align*} (\id\tens&\underline\Delta)\Theta=S(\CR^2\CR^{'1})\tens \CR^1{}_1\CR^{'2}{}_1S\CR^{iv2}\tens\CR^{iv1}{}_1(\CR^1{}_2\CR^{'2}{}_2)S\CR^{iv}{}_2\\
&=S(\CR^2\CR^{''2}\CR^{'''1}\CR^{'1})\tens\CR^1\CR^{'2}S(\CR^{iv2}\CR^{v2})\tens \CR^{iv1}\CR^{''1}\CR^{'''2}S\CR^{v1}\\
&=S(\CR^2\CR^{''2}\CR^{'''1}\CR^{'1})\tens\CR^1\CR^{'2}\CR^{v2}S(\CR^{iv2})\tens \CR^{iv1}\CR^{''1}\CR^{'''2}\CR^{v1}\\
&=S(\CR^2\CR^{''2}\CR^{'1}\CR^{'''1})\tens\CR^1\CR^{v2}\CR^{'2}S(\CR^{iv2})\tens \CR^{iv1}\CR^{''1}\CR^{v1}\CR^{'''2}\\
&=S(\CR^2\CR^{'1}\CR^{''2}\CR^{'''1})\tens\CR^1\CR^{'2}\CR^{v2}S(\CR^{iv2})\tens \CR^{iv1}\CR^{v1}\CR^{''1}\CR^{'''2}\\
&=S(\CR^{''2}\CR^{'''1})S(\CR^2\CR^{'1})\tens\CR^1\CR^{'2}\tens\CR^{''1}\CR^{'''2}=(\mu\tens\id)\Theta_{12}\Theta_{23}
\end{align*} 
where for the 3rd equality we used that $(S\tens S)(\CR)=\CR$. The other steps are similar to the previous case. \endproof

For the associated Hadamard gate we recover $\ch^{-1}=\CS$  the braided Fourier transform $\underline H\to \underline H$ in \cite{LyuMa}, where it was shown in the case of a ribbon Hopf algebra, such as $u_q(sl_2)$,  that $\CS$ along with left multiplication $\CT$ by the ribbon element, obey the modular identity
\begin{equation}\label{mod} (\CS\CT)^3=\lambda\CS^2\end{equation}
for some constant $\lambda$. This is key to the construction of topological invariants from quantum groups. One also has $\CS^2=\underline S^{-1}$ and hence
\begin{equation}\label{modh2} \mathfrak h^2=\underline S\end{equation}
is the braided antipode. Finally, in the case where $H$ is a quasitriangular flip-Hopf $*$-algebra with $\CR^\dagger=\CR^{-1}$, which is the case for $u_q(sl_2)$ in Example~\ref{exsl2}, the category of $H$-modules becomes a braided bar category, see  \cite[Prop 3.7, Thm~3.6]{BegMa:bar} for details. In the ribbon case, this is a strong one (meaning that the natural equivalence between applying the bar functor twice and the identity respects tensor products). It is then proven \cite[Thm~5.7]{BegMa:bar} that $\underline H$ becomes a braided $*$-Hopf algebra in the category with a modified operation 
\begin{equation}\label{transtar}\underline{*}:\underline H\to \overline{\underline H},\quad  h^{\underline *}=(S^2\CR^1) h^* \CR^2v;\quad v=\CR^1 S\CR^2,\end{equation}
where the element $v\in H$ implements $S^{-2}$ by conjugation as part of Drinfeld theory\cite{Dri} and the output of $\underline{*}$ is viewed in the conjugate object built on the same vector space. Note that in this bar category of $H$-modules, the natural equivalence $\Upsilon:\overline{ V\tens W}\to \overline W\tens\overline{V}$ is, like the braiding,  governed by $\CR$ and is not the usual flip map as for vector spaces, hence $\underline{*}$ is not a usual $*$-involution. We do however, have $\underline{*}\circ\underline S=\underline S^{-1}\circ\underline{*}$, see \cite[Prop.~5.2]{BegMa:bar}. We limit ourselves to one new observation.
\begin{lemma} For $H$ a quasitriangular flip-Hopf $*$-algebra with $\CR^\dagger=\CR^{-1}$, we have
\[ \underline S\circ \underline{*}= S\circ *\]
i.e. the underlying antilinear automorphism $\theta=S\circ*$ on $H$ is unchanged by transmutation. 
\end{lemma}
\proof From the definitions,
\begin{align*}\underline S(h^{\underline*})&=\CR^2S(\CR^1{}_1(S^2\CR^{''1})h^*\CR^{''2}vS\CR^1{}_2=\CR^2\CR^{'2}S(\CR^1 (S^2\CR^{''1})h^*\CR^{''2}v S\CR^{'1})\\
&=\CR^2(\CR^{'2}S^2\CR^{'1})(Sv)\CR^{''2}(Sh^*)S(\CR^1S\CR^{''1})=\CR^2\CR^{''2}(Sh^*)S(\CR^1S\CR^{''1})=Sh^*\end{align*}
where $\CR^2S\CR^1=u^{-1}$ in Drinfeld theory and $u=S v$, see \cite[Prop.~2.1.8]{Ma}.\endproof

This suggests that a braided version of Section~\ref{secstar} should be possible, even though significantly more complicated due to working in a nontrivial bar category. 

\section{Concluding remarks}

An open question at the algebra level is whether the notion of F-Hopf algebra or its braided version go truly beyond those obtained by amplifying a single Hopf algebra as in Propositions~\ref{propFhopf},~\ref{propbraFhopf}. In principle, the red and green Frobenius structures need not be closely related, with the result that the two Hopf algebras need not be op-dual to each other, but no examples are known even in the unbraided case. We also saw that the simplest version of the $*$-algebra theory in Section~\ref{secuni}, where $(\ ,\ )$ combines with a $*$-structure to give a Hilbert space, does not apply to non-unimodular quantum groups such as $u_q(sl_2)$ with $q$ a primitive odd root of unity. This means that the natural sesquilinear form defined by $\int h^*h'$ does not make $H$ into a Hilbert space since $\overline{\int h^*h'}=\int S( h'{}^*h)$ is not necessarily $\int h'{}^*h$. We saw at the end of Section~\ref{sectrans} that these problems are even more pronounced in the braided case where the natural transmuted $\underline{*}$ on $\underline{H}$, while it makes the latter into a $*$-Hopf algebra in the relevant braided bar category, is not a $*$-algebra structure in the usual sense. This topic of $*$-structures merits further investigation, perhaps guided by is its precise role in quantum computing,  which is currently unclear. In particular, we have not developed a general theory for when the `gates' in question are unitary, which would be needed in applications. The situation is rather better for a Drinfeld double $H=D(G)$ of a finite group $G$, which is unimodular. 

Another feature of the paper was an abstract notion of `Hadamard gate'. Its inverse appears (\ref{hadfou}) as Hopf algebra Fourier transform composed evaluated against an element $\Theta\in H\tens H$ give an operator. We identified three types, with Type 1 connecting red and green Frobenius structures, which fits well with their usual role in  ZX calculus, and Type 2 (resp. Type 3) corresponding to Hopf algebra (anti) self duality of $H$, which fits better with the algebraic picture (where the map on coproducts would naturally go the other way to the map on products related by duality). The practicality of this in quantum computing would need to be looked at further. The theory also applies in the braided case if we have integrals. Here, braided Fourier transform is also of interest in the noncommutative geometry of quantum groups\cite{Ma:hod}, but one has to deal with the fact that the integral is often not a morphism to the trivial object, even on some very basic examples such as the braided line. It would be interesting to adapt the approach to the braided line in \cite{Ma:hod} so as to still obtain a braided F-Hopf algebra, which we saw did not otherwise work. We saw that there is no problem, however, for the braided theory in the case of braided Hopf algebras obtained by transmutation from a factorisable quantum group such as $u_q(sl_2)$ at odd roots of unity. It is also striking that the braided self-duality/Type 2 Hadamard gate that applies canonically in this case obeys\cite{LyuMa} the modular identities in (\ref{mod}), making it central to one approach to quantum group topological invariants. 

Finally, while the present constructions are interesting from an algebraic point of view, it remains to be seen how useful they are in quantum computing. The general idea is that while standard ZX-calculus is well adapted to conventional quantum computers, these suffer in practice from noise problems resulting in a current focus is on the development of `topologically fault tolerant' methods. The prototype of these is the Kitaev model\cite{Kit} for which the underlying algebraic structure is a  Drinfeld double $D(G)$. This can be generalised to the quantum double of a finite-dimensional Hopf $C^*$-algebras\cite{BMCA} and beyond\cite{CowMa}. As noted in \cite{Meu}, such models are closer in spirit to the Turaev-Viro invariants of 3-manifolds,  which have underlying them the quantum double of $u_q(sl_2)$, while the related Jones/Reshetkhin-Turaev knot invariant come from $u_q(su_2)$.  The fact that modular identities (\ref{mod}) for the Turaev-Viro invariant have a natural role as Hadamard gates loosely supports the possibility of a braided ZX-calculus approach to topological quantum computing. This could be explored further. Here, Fibonnaci anyons of interest in fault tolerant quantum computing are related to $u_q(sl_2)$ at a 5th root of unity\cite{fib}.

Another setting in which all of these ideas could potentially be explored concretely is `digital quantum computing' as proposed in \cite{MaPac}. This work classifies Hopf algebras over $\F_2$ in low dimension, including which of them are quasitriangular. One could then search for digital Hadamard forms as well look at digital braided versions.

\end{document}